\documentclass[nofonts]{tufte-handout}
\usepackage[utf8]{inputenc}
\usepackage{appendix}
\usepackage{import}
\usepackage{hyperref}

\usepackage{import}
\usepackage{blindtext}

\usepackage{graphicx}
\setkeys{Gin}{width=\linewidth,totalheight=\textheight,keepaspectratio}
\usepackage[caption=false]{subfig}

\usepackage{fancyhdr}
\usepackage{amsmath,amssymb,amsthm}
\usepackage{tikz}
\usepackage{wrapfig} 
\usetikzlibrary{cd, arrows,automata,positioning}
\usepackage{tipa} 
\usepackage{mathtools} 
\usepackage{caption}
\usepackage{thmtools, thm-restate} 

\usetikzlibrary{decorations.pathreplacing} 

\DeclareMathOperator{\Aa}{\mathcal{A}}

\DeclareMathOperator{\Ca}{\mathcal{C}}

\DeclareMathOperator{\Oa}{\mathcal{O}}

\DeclareMathOperator{\Pb}{\mathbb{P}}

\DeclareMathOperator{\Rb}{\mathbb{R}}

\DeclareMathOperator{\Zb}{\mathbb{Z}}

\DeclareMathOperator{\iso}{\simeq}
\DeclareMathOperator{\from}{\leftarrow}

\DeclareMathOperator{\op}{^\text{op}}
\DeclareMathOperator{\inv}{^{-1}}
\DeclareMathOperator{\id}{\mathsf{id}}

\DeclareMathOperator{\ob}{\mathsf{ob}}

\DeclareMathOperator{\Set}{\mathbf{Set}}

\DeclareMathOperator{\Cat}{\mathbf{Cat}}

\DeclareMathOperator{\Finset}{\mathbf{FinSet}}
\DeclareMathOperator{\Dynam}{\mathsf{Dynam}}
\DeclareMathOperator{\Cospan}{\mathsf{Cospan}}
\DeclareMathOperator{\Span}{\mathsf{Span}}

\DeclareMathOperator{\Mnfld}{\mathsf{Mnfld}}
\DeclareMathOperator{\Lens}{\mathsf{Lens}}

\usepackage{stmaryrd}
\newcommand{\cp}{\mathbin{\fatsemi}}

\newcommand{\xto}[1]{\xrightarrow{#1}}

\newcommand{\pto}{\,\cdot\kern-.1em{\to}\,}

\makeatletter
\providecommand*{\xmapstofill@}{%
  \arrowfill@{\mapstochar\relbar}\relbar\rightarrow
}
\providecommand*{\xmapsto}[2][]{%
  \ext@arrow 0395\xmapstofill@{#1}{#2}%
}

\makeatletter
\def\slashedarrowfill@#1#2#3#4#5{%
  $\m@th\thickmuskip0mu\medmuskip\thickmuskip\thinmuskip\thickmuskip
   \relax#5#1\mkern-7mu%
   \cleaders\hbox{$#5\mkern-2mu#2\mkern-2mu$}\hfill
   \mathclap{#3}\mathclap{#2}%
   \cleaders\hbox{$#5\mkern-2mu#2\mkern-2mu$}\hfill
   \mkern-7mu#4$%
}
\def\rightslashedarrowfill@{%
  \slashedarrowfill@\relbar\relbar\mapstochar\rightarrow}
\newcommand\xslashedrightarrow[2][]{%
  \ext@arrow 0055{\rightslashedarrowfill@}{#1}{#2}}
\makeatother

\newcommand{\topro}{\xslashedrightarrow{}}

\tikzset{
    vert/.style={anchor=south, rotate=90, inner sep=.5mm}
} 

\newtheorem{thm}{Theorem}[section]

\theoremstyle{definition}
\newtheorem{defn}[thm]{Definition}

\newtheorem{ex}{Example}

\newtheorem{lem}[thm]{Lemma}

\DeclareMathOperator{\Wd}{\mathsf{WD}}

\DeclareMathOperator{\CDS}{\mathsf{CDS}}

\DeclareMathOperator{\Rsm}{\mathsf{RSM}}

\newcommand{\lens}[2]{{{#1} \choose {#2}}}
\newcommand{\lensbig}[2]{\begin{pmatrix} {#1} \\ {#2} \end{pmatrix}}

\newcommand{\lensmap}[2]{\overset{#1}{\underset{#2}{\leftrightarrows}}}

\newcommand{\biKl}{\mathsf{BiKleisli}}
\newcommand{\Sub}{\mathsf{Sub}}
\newcommand{\Riem}{\mathsf{Riem}}
\newcommand{\Euc}{\mathsf{Euc}}

\usepackage{amsthm}
\usepackage{stmaryrd}
\usepackage{newpxtext}
\usepackage[varg,bigdelims]{newpxmath}
\usepackage[usenames,dvipsnames]{xcolor}
\usepackage{tikz}

\usetikzlibrary{
	cd,
	math,
	backgrounds,
	decorations.markings,
	decorations.pathreplacing,
	positioning,
	arrows.meta,
	circuits.logic.US,
	shapes,
	calc,
	fit,
	trees,
	quotes}

\tikzset{
  WD/.style={
  	label/.style={
    	font=\everymath\expandafter{\the\everymath\scriptstyle},
      inner sep=0pt,
      node distance=2pt and -2pt},
  	label distance=-2pt,
  	every to/.style={draw},
    semithick,
    node distance=\bbx and \bby,
    decoration={markings, mark=at position \stringdecpos with \stringdec},
    bb port length=3pt,
  	bb port sep=1,
		bb inside color=white,
		bb outside color=black,
	 	bbx = .4cm,
		bb min width=.4cm,
	  bby = 2ex,
	  bb penetrate=0,
	  bb rounded corners=2pt,
	  dot size=3pt,
    shell size = 16pt,
   	penetration = 0pt,
    link size = 2pt,
    shell color = blue,
  	shell inside color=\pcolor!20,
 	  shell outside color=\pcolor!50!black,
  	surround sep=2pt,
    ar/.style={postaction={decorate}},
  	execute at begin picture={\tikzset{
  		x=\bbx, y=\bby, 
			circuit logic US, tiny circuit symbols
			}
		}
  },
  beamer/.style={
  	bbx=.4cm,
		bb min width=.4cm,
		bby=6pt,
		bb port length=3pt,
		bb port sep=.5,
  	dot size=1pt,
    shell size = 11pt, 
   	penetration = 0pt,
    link size = 1pt,
    shell color = blue,
    surround sep=1pt,
    inner sep=1pt,
    font=\tiny,
    bb inside color=\picolor,
    bb outside color=\pocolor,
	},
	bb standard colors/.style={bb inside color=white, bb outside color=black},
	bb inside color/.store in=\bbicolor,
	bb outside color/.store in=\bbocolor,
  bbx/.store in=\bbx,
  bby/.store in=\bby,
  bb port sep/.store in=\bbportsep,
  bb port length/.store in=\bbportlen,
  bb penetrate/.store in=\bbpenetrate,
  bb min width/.store in=\bbminwidth,
  bb rounded corners/.store in=\bbcorners,
  bb/.code 2 args={
    \pgfmathsetlengthmacro{\bbheight}{\bbportsep * (max(#1,#2)+1) * \bby}
    \pgfkeysalso{
      draw=\bbocolor,
      fill=\bbicolor,
      minimum height=\bbheight,
      minimum width=\bbminwidth,
      outer sep=0pt,
      rounded corners=\bbcorners,
      thick,
      prefix after command={\pgfextra{\let\fixname\tikzlastnode}},
      append after command={\pgfextra{\draw
      	\ifnum #1=0{} \else foreach \i in {1,...,#1} {
        	($(\fixname.north west)!{(2*\i-1)/(2*#1)}!(\fixname.south west)$) +(-\bbportlen,0) coordinate (\fixname_in\i) -- +(\bbpenetrate,0) coordinate (\fixname_in\i')}\fi 
        \ifnum #2=0{} \else foreach \i in {1,...,#2} {
        	($(\fixname.north east)!{(2*\i-1)/(2*#2)}!(\fixname.south east)$) +(-
\bbpenetrate,0) coordinate (\fixname_out\i') -- +(\bbportlen,0) coordinate (\fixname_out\i)}\fi
;
       }}}
		},
	dot size/.store in=\dotsize,
	dot/.style={
		circle, draw, thick, inner sep=0, fill=black, minimum width=\dotsize
	},
	bb name/.style={
    append after command={
		\pgfextra{\node[anchor=north] at (\fixname.north) {#1};}
		}
	},
  shell size/.store in=\psize,
	penetration/.store in=\penetration,
  spacing/.store in=\spacing,
  link size/.store in=\lsize,
  shell color/.store in=\pcolor,
 	shell inside color/.store in=\picolor,
 	shell outside color/.store in=\pocolor,
 	surround sep/.store in=\ssep,
 	link/.style={
  	circle, 
  	draw=black, 
  	fill=black,
  	inner sep=0pt, 
 		minimum size=\lsize
 	},
  shell/.style={
 		circle, 
 		draw = \pocolor, 
  	fill = \picolor,
  	minimum size = \psize
  },
  func/.style={
  	shell,
		rectangle,
		rounded corners=.5*\psize,
		inner ysep=.125*\psize,
		minimum width=1.125*\psize,
		inner xsep=.25*\psize,
  },
  funcr/.style={
    func,
    rectangle round north west=false, 
		rectangle round south west=false,
  },
  funcl/.style={
    func,
		rectangle round north east=false, 
		rectangle round south east=false,
  },
  funcu/.style={
    func,
		rectangle round south east=false, 
		rectangle round south west=false,
  },
  funcd/.style={
    func,
		rectangle round north east=false, 
		rectangle round north west=false,
  },
  outer shell/.style={
 		ellipse, 
 		draw,
  	inner sep=\ssep,
  	color=gray,
 	},
  intermediate shell/.style={
 		ellipse,
 		dashed, 
  	draw,
  	inner sep=\ssep,
 		color=\pocolor,
 	},
 }

\tikzset{
	oriented WD/.style={
		every to/.style={out=0,in=180,draw},
    label/.style={
    	font=\everymath\expandafter{\the\everymath\scriptstyle},
      inner sep=0pt,
      node distance=2pt and -2pt},
    semithick,
    node distance=1 and 1,
    decoration={markings, mark=at position \stringdecpos with \stringdec},
    ar/.style={postaction={decorate}},
    execute at begin picture={\tikzset{
    	x=\bbx, y=\bby,
      every fit/.style={inner xsep=\bbx, inner ysep=\bby}}}
    },
    string decoration/.store in=\stringdec,
    string decoration={\arrow{stealth};},
    string decoration pos/.store in=\stringdecpos,
    string decoration pos=.7,
    bbx/.store in=\bbx,
    bbx = 1.5cm,
    bby/.store in=\bby,
    bby = 1.ex,
    bb port sep/.store in=\bbportsep,
    bb port sep=1.5,
    bb port length/.store in=\bbportlen,
    bb port length=4pt,
    bb penetrate/.store in=\bbpenetrate,
    bb penetrate=0,
    bb min width/.store in=\bbminwidth,
    bb min width=1cm,
    bb min height/.store in=\bbminheight,
    bb min height=1cm,
    bb rounded corners/.store in=\bbcorners,
    bb rounded corners=2pt,
    bb spider/.style={
    	bb port sep=1, bb port length=10pt, bbx=.4cm, bb min width=.4cm, bby=.8ex},
    bb small/.style={
    	bb port sep=1, bb port length=2.5pt, bbx=.4cm, bb min width=.4cm, bby=.7ex},
		bb medium/.style={
			bb port sep=1, bb port length=2.5pt, bbx=.4cm, bb min width=.4cm, bby=.9ex},
    bb/.code n args={4}{
    	\pgfmathsetlengthmacro{\bbheight}{\bbportsep * (max(#1,#2)+1) * \bby}
    	\pgfmathsetlengthmacro{\bbwidth}{\bbportsep * (max(#3,#4)+1) * \bby}
      \pgfkeysalso{draw,minimum height=\bbminheight,minimum
       width=\bbminwidth,outer sep=0pt,
         rounded corners=\bbcorners,thick,
         prefix after command={\pgfextra{\let\fixname\tikzlastnode}},
         append after command={\pgfextra{\draw
            \ifnum #1=0{} \else foreach \i in {1,...,#1} {
            	($(\fixname.north west)!{\i/(#1+1)}!(\fixname.south west)$) +(-\bbportlen,0) coordinate (\fixname_in\i) -- +(\bbpenetrate, 0) coordinate (\fixname_in\i')}\fi 
            \ifnum #2=0{} \else foreach \i in {1,...,#2} {
            	($(\fixname.north east)!{\i/(#2+1)}!(\fixname.south east)$) +(-
\bbpenetrate,0) coordinate (\fixname_out\i') -- +(\bbportlen,0) coordinate (\fixname_out\i)}\fi

\ifnum #3=0{} \else foreach \i in {1,...,#3} {
            	($(\fixname.north west)!{\i/(#3+1)}!(\fixname.north east)$) +(0, \bbportlen) coordinate (\fixname_top\i) -- +(0,-\bbpenetrate) coordinate (\fixname_top\i')}\fi 
\ifnum #4=0{} \else foreach \i in {1,...,#4} {
            	($(\fixname.south west)!{\i/(#4+1)}!(\fixname.south east)$) +(0, \bbpenetrate) coordinate (\fixname_bot\i) -- +(0,-\bbportlen) coordinate (\fixname_bot\i')}\fi 
;
           }}}
		},
			bb name/.style={
     	append after command={
				\pgfextra{\node[anchor=north] at (\fixname.north) {#1};}
			}
		}
  }

  \tikzset{
  	unoriented WD/.style={
  		every to/.style={draw},
  		shorten <=-\penetration, shorten >=-\penetration,
  		label distance=-2pt,
  		thick,
  		node distance=\spacing,
  		execute at begin picture={\tikzset{
  			x=\spacing, y=\spacing}}
  		},
  	pack size/.store in=\psize,
  	pack size = 8pt,
  	spacing/.store in=\spacing,
  	spacing = 8pt,
  	link size/.store in=\lsize,
  	link size = 2pt,
		penetration/.store in=\penetration,
		penetration = 2pt,
  	pack color/.store in=\pcolor,
  	pack color = blue,
  	pack inside color/.store in=\picolor,
  	pack inside color=blue!20,
  	pack outside color/.store in=\pocolor,
  	pack outside color=blue!50!black,
  	surround sep/.store in=\ssep,
  	surround sep=8pt,
  	link/.style={
  		circle, 
  		draw=black, 
  		fill=black,
  		inner sep=0pt, 
  		minimum size=\lsize
  	},
  	pack/.style={
  		circle, 
  		draw = \pocolor, 
  		fill = \picolor,
  		inner sep = .25*\psize,
  		minimum size = \psize
  	},
  	outer pack/.style={
  		ellipse, 
  		draw,
  		inner sep=\ssep,
  		color=\pocolor,
  	},
  	intermediate pack/.style={
  		ellipse,
  		dashed, 
  		draw,
  		inner sep=\ssep,
  		color=\pocolor,
  	},
  }

\tikzset{
	spider diagram/.style={
		every to/.style={out=0, in=180, draw, thick},
		thick
	},
	dot size/.store in=\dotsize,
	dot size = 5pt,
	dot fill/.store in=\dotfill,
	dot fill = black,
	leg length/.store in=\leglen,
	leg length = 15pt,
	baby/.style={dot size = 2pt, leg length = 6pt},
	young/.style={dot size = 3pt, leg length = 10pt},
	special spider/.code n args={4}{
		\pgfkeysalso{circle, draw, thick, inner sep=0, fill=\dotfill, minimum width=\dotsize,
  		prefix after command={\pgfextra{\let\fixname\tikzlastnode}},
  		append after command={\pgfextra{
  			\ifnum #1=0{} \else {\foreach \i in {1,...,#1} {
					\tikzmath{\anglei={-90*(#1+1-2*\i)/#1};}
  				\draw [thick]
						(\fixname) .. controls 
						($(\fixname.center)-(\anglei:#3/3)$) and ($(\fixname.center)-(\anglei:#3*2/3)$) .. 
						({$(\fixname)-(\anglei:#3*2/3)$}-|{$(\fixname)-(#3,0)$}) coordinate (\fixname_in\i);
  			}}\fi
  			\ifnum #2=0{} \else {\foreach \i in {1,...,#2} {
					\tikzmath{\anglei={90*(#2+1-2*\i)/#2};}
  				\draw [thick]
						(\fixname.center) .. controls 
						($(\fixname.center)+(\anglei:#4/3)$) and ($(\fixname.center)+(\anglei:#4*2/3)$) .. 
						({$(\fixname.center)+(\anglei:#4*2/3)$}-|{$(\fixname.center)+(#4,0)$}) coordinate (\fixname_out\i);
  			}}\fi
  		}}
		}
	},
	spider/.code 2 args={
		\pgfkeysalso{special spider={#1}{#2}{\leglen}{\leglen}}
	}
}

\tikzset{Yonepart/.pic={
	\node[bb={1}{2},bb name = {\tiny$X_{11}$}] (X11) {};
	\node[bb={2}{2},below right=of X11,bb name = {\tiny$X_{12}$}] (X12) {};
	\node[bb={2}{1}, above right=of X12,bb name = {\tiny$X_{13}$}] (X13) {};
	\node[bb={2}{2}, fit={($(X11.north west)+(.3,1.5)$) (X12)  ($(X13.east)+(-.3,0)$)},bb name = {\scriptsize $Y_1$}] (Y1) {};
	\draw (Y1_in1') to (X11_in1);	
	\draw (Y1_in2') to (X12_in2);
	\draw (X11_out1) to (X13_in1);
	\draw (X11_out2) to (X12_in1);
	\draw (X12_out1) to (X13_in2);
	\draw (X12_out2) to (Y1_out2');
	\draw (X13_out1) to (Y1_out1');
	\coordinate (bottombox) at ($(X12.south)$);
	\coordinate (rightbox) at ($(X13.east)$);
	\coordinate (Y1northwest) at ($(Y1.north west)$);
	}
}

\tikzset{Ytwopart/.pic={
	\node[bb={2}{2}, bb name = {\tiny$X_{21}$}] (X21) {};
	\node[bb={1}{2},above right=-1 and 1 of X21,bb name = {\tiny$X_{22}$}] (X22) {};
	\node[bb={1}{2}, fit={($(X21.south west)+(-.25,0)$) ($(X22.north east)+(.25,3.5)$)},bb name = {\scriptsize$Y_2$}] (Y2){};
	\draw (Y2_in1') to (X21_in2);
	\draw (X21_out1) to (X22_in1);
	\draw (X22_out2) to (Y2_out1');
	\draw let \p1=(X22.south east), \p2=($(Y2_out2)$), \n1={\y1-\bby}, \n2=\bbportlen in
	  (X21_out2) to (\x1+\n2,\n1) -- (\x1+\n2,\n1) to (Y2_out2');
	\draw let \p1=(X22.north east), \p2=(X21.north west), \n1={\y1+\bby}, \n2=\bbportlen in
          (X22_out1) to[in=0] (\x1+\n2,\n1) -- (\x2-\n2,\n1) to[out=180] (X21_in1);
          }
}

\tikzset{SmallNeuronPic/.pic={
 \node[bb={3}{1}] (N1) {$\scriptstyle N_1$};
  \node[bb={2}{1}, above right=.7 and 3.5 of N1] (N2) {$\scriptstyle N_2$};
  \node[bb={2}{1}, below =of N2] (N3) {$\scriptstyle N_3$};
  \node[bb={3}{1}, below =of N3] (N4) {$\scriptstyle N_4$};
  \node[bb={6}{8}, fit={($(N1.west)-(.5,0)$) ($(N2.north)+(0,2)$) ($(N3.east)+(1.5,0)$) ($(N4.south)-(0,1)$)}, bb name={$\scriptstyle X$}] (X) {};
  \draw (X_in1') to (N2_in1);
  \draw (X_in2') to (N1_in1);
  \draw (X_in3') to (N1_in2);
  \draw (X_in4') to (N1_in3);
  \draw (X_in6') to (N4_in2);
  \draw (N1_out1) to (N2_in2);
  \draw (N1_out1) to (N3_in1);
  \draw (N1_out1) to (N4_in1);
  \draw (N2_out1) to (X_out1');
  \draw (N2_out1) to (X_out2');
  \draw (N2_out1) to (X_out3');
  \draw (N3_out1) to (X_out4');
  \draw (N3_out1) to (X_out5');
  \draw (N3_out1) to (X_out6');
  \draw (N4_out1) to (X_out7');
  \draw (N4_out1) to (X_out8'); 
  \draw (X_in5') to[looseness=2] (N3_in2);
  \draw let \p1=(N4.south east), \p2=(N4.south west), \n1={\y2-\bby}, \n2=\bbportlen in
          (N3_out1) to[in=0] (\x1+\n2,\n1) -- (\x2-\n2,\n1) to[out=180] (N4_in3);
}
}

\tikzset{SmallNeuronDashed/.pic={
 \node[bb={3}{1}] (N1) {$\scriptstyle N_1$};
  \node[bb={2}{1}, above right=.7 and 3.5 of N1] (N2) {$\scriptstyle N_2$};
  \node[bb={2}{1}, below =of N2] (N3) {$\scriptstyle N_3$};
  \node[bb={3}{1}, below =of N3] (N4) {$\scriptstyle N_4$};
  \node[bb={6}{8}, fit={($(N1.west)-(.5,0)$) ($(N2.north)+(0,2)$) ($(N3.east)+(1.5,0)$) ($(N4.south)-(0,1)$)}, bb name={$\scriptstyle X$}] (X) {};
  \draw[dashed] (X_in1') to (N2_in1);
  \draw[dashed] (X_in2') to (N1_in1);
  \draw[dashed] (X_in3') to (N1_in2);
  \draw[dashed] (X_in4') to (N1_in3);
  \draw[dashed] (X_in6') to (N4_in2);
  \draw[dashed] (N1_out1) to (N2_in2);
  \draw[dashed] (N1_out1) to (N3_in1);
  \draw[dashed] (N1_out1) to (N4_in1);
  \draw[dashed] (N2_out1) to (X_out1');
  \draw[dashed] (N2_out1) to (X_out2');
  \draw[dashed] (N2_out1) to (X_out3');
  \draw[dashed] (N3_out1) to (X_out4');
  \draw[dashed] (N3_out1) to (X_out5');
  \draw[dashed] (N3_out1) to (X_out6');
  \draw[dashed] (N4_out1) to (X_out7');
  \draw[dashed] (N4_out1) to (X_out8'); 
  \draw[dashed] (X_in5') to[looseness=2] (N3_in2);
  \draw[dashed] let \p1=(N4.south east), \p2=(N4.south west), \n1={\y2-\bby}, \n2=\bbportlen in
          (N3_out1) to[in=0] (\x1+\n2,\n1) -- (\x2-\n2,\n1) to[out=180] (N4_in3);
}
}

\tikzset{SmallNestingPic/.pic={
\path (0,0) pic [purple] {Yonepart};
\path ($(rightbox)+(5,-5)$) pic [orange] {Ytwopart};
 
\node[bb={1}{2}, fit={($(Y1northwest)+(-.5,4)$) ($(Y2.south east)+(1,0)$)}, bb name={\small $Z$}] (Z) {};
\draw (Z_in1') to (Y1_in2);
\draw let \p1=(Y2.north west),\p2=(Y2.north east),\n1={\y2+\bby},\n2=\bbportlen in
  (Y1_out1) to (\x1+\n2,\n1)--(\x2+\n2,\n1) to (Z_out1');
\draw (Y1_out2) to (Y2_in1);
\draw (Y2_out2) to (Z_out2');
\draw let \p1=(Y2.north east), \p2=(Y1.north west), \n1={\y2+\bby}, \n2=\bbportlen in
          (Y2_out1) to[in=0] (\x1+\n2,\n1) -- (\x2-\n2,\n1) to[out=180] (Y1_in1);
          }
}

\tikzset{Zredgreen/.pic={
\node[bb={2}{2}, green!50!black, bb name = $\scriptstyle Y_1$] (YY1) {};
\node[bb={1}{2}, red, below right=-1 and 2 of YY1, bb name=$\scriptstyle Y_2$] (YY2) {};
\node[bb={1}{2}, fit={($(YY1.north west)+(-.5,4)$) ($(YY2.south east)+(.5,-2)$)}, bb name={\scriptsize $Z$}] (Z) {};
\draw (Z_in1') to (YY1_in2);
\draw (YY1_out1) to (Z_out1');
\draw (YY1_out2) to (YY2_in1);
\draw (YY2_out2) to (Z_out2');
\draw let \p1=(YY2.north east), \p2=(YY1.north west), \n1={\y2+\bby}, \n2=\bbportlen in
          (YY2_out1) to[in=0] (\x1+\n2,\n1) -- (\x2-\n2,\n1) to[out=180] (YY1_in1);
}
}

\tikzset{Zcombined/.pic={
	\node[bb={1}{2},green!25!black,bb name = {\tiny$X_{11}$}] (X11) {};
	\node[bb={2}{2},green!25!black,below right=of X11,bb name = {\tiny$X_{12}$}] (X12) {};
	\node[bb={2}{1}, green!25!black,above right=of X12,bb name = {\tiny$X_{13}$}] (X13) {};
	\draw (X11_out1) to (X13_in1);
	\draw (X11_out2) to (X12_in1);
	\draw (X12_out1) to (X13_in2);

	\node[bb={2}{2}, red!30!black, below right = 0 and 1.25 of X12, bb name = {\tiny$X_{21}$}] (X21) {};
	\node[bb={1}{2}, red!30!black, above right=-1 and 1 of X21,bb name = {\tiny$X_{22}$}] (X22) {};
	\draw (X21_out1) to (X22_in1);
	\draw let \p1=(X22.north east), \p2=(X21.north west), \n1={\y1+\bby}, \n2=\bbportlen in
          (X22_out1) to[in=0] (\x1+\n2,\n1) -- (\x2-\n2,\n1) to[out=180] (X21_in1);
        
        \node[bb={1}{2}, fit = {($(X11.north east)+(-1,3)$) (X12) (X13) ($(X21.south)+(0,-1)$) ($(X22.east)+(.5,0)$)}, bb name ={\scriptsize $Z$}] (Z) {};
	
	\draw (Z_in1') to (X12_in2);
	\draw (X13_out1) to (Z_out1');
	\draw (X12_out2) to (X21_in2);
	\draw let \p1=(X22.south east),\n1={\y1-\bby}, \n2=\bbportlen in
	  (X21_out2) to (\x1+\n2,\n1) to (Z_out2');
	\draw let \p1=(X22.north east), \p2=(X11.north west), \n1={\y2+\bby}, \n2=\bbportlen in
          (X22_out2) to[in=0] (\x1+\n2,\n1) -- (\x2-\n2,\n1) to[out=180] (X11_in1);
}
}

\title[Resource Sharing Machines]{{\normalfont \huge An Algebra of Resource Sharing Machines}\\ {\normalfont \LARGE Unifying Two Flavors of Open Dynamical Systems}}

\author{Sophie Libkind}

\begin{document}

\noindent
\makebox[\textwidth][l]{%
  \begin{minipage}[b]{\textwidth+\marginparsep+\marginparwidth}
  \raggedright
  \maketitle
  \end{minipage}%
}


\begin{abstract}
\noindent {\bf Abstract} Dynamical systems are a broad class of  mathematical tools used to describe the evolution of physical and computational processes. Traditionally these processes model changing entities in a static world. Picture a ball rolling on an empty table. In contrast, \emph{open} dynamical systems model changing entities in a \emph{changing} world. Picture a ball in an ongoing game of billiards. In the literature, there is ambiguity about the interpretation of the "open" in open dynamical systems. In other words, there is ambiguity in the mechanism by which open dynamical systems interact. To some, open dynamical systems are input-output machines which interact by feeding the input of one system with the output of another. To others, open dynamical systems are input-output agnostic and interact through a shared pool of resources. 

In this paper, we define an algebra of open dynamical systems which unifies these two perspectives. We consider in detail two concrete instances of dynamical systems --- continuous flows on manifolds and non-deterministic automata.
\end{abstract}

\section{Introduction}

Classical information and communication theory assumes a one way communication channel between an established sender and receiver. In contrast, physical systems do not have a mechanism for such directed interaction --- a property  caricatured by the slogan "every action has an equal and opposite reaction." These observations lead to the following mystery: how does a physical system (such as a transistor or cell)  reliably represent an ideal computer (such as a NOT gate or  gene regulatory network)? As a first step towards solving this puzzle, we give a general framework for system interaction which captures both undirected and uni-directional communication. These distinct types of interactions are captured by mathematical notions we respectively call \emph{resource sharers}~\citep{ baez2016compositional, Baez_2017}  and \emph{machines}~\citep{vagner2014algebras, schultz2016dynamical, spivak2020poly, myers2020double}.

Dynamical systems refer to a broad class mathematical objects which model "things that change." A Turing machine (and more generally, a computer) is a dynamical system; the state of the tape changes according to an algorithm. An electromagnetic field is also a dynamical system; the state of the field changes  according to the laws of Gauss, Faraday, and Maxwell. Traditionally, mathematicians and scientists fix a dynamical system and then ask questions about it. What are its equilibrium points? Its orbits? Its entropy? However, in nature dynamical systems do not exist in isolation. For example, the state of a computer is influenced by the network of servers it is connected to and the actions of the user on the keyboard and mouse.  Hence, we are interested in \emph{open} dynamical systems, i.e. those which have a mechanism for interacting with other systems. When dynamical systems interact we say they \emph{compose}. Resource sharers and machines are two flavors of open dynamical systems with distinct styles of composition.

Inspired by the physical interactions in chemical reaction networks, the authors of~\citep{Baez_2017} define a framework for open dynamical systems which compose as resource sharers. Two resource sharers compose by simultaneously affecting and reacting to a shared pool of resources. When resource sharers compose:
\begin{enumerate}
    \item \emph{communication is undirected}. Each system may both affect and be affected by the state of the shared pool of resources. Through this medium, they may both affect and be affected by each other.
    \item \emph{interaction is passive}. The communication channel is incidental to the  fact that the systems refer to the same resource. The rules for how each system affects and reacts to state of the resource is independent of the action of other systems on the pool. 
\end{enumerate}
When people communicate verbally, they are composing as resource sharers where the shared resource is "vibrations in the air space." All participants in a conversation affect and are affected by the changing state of the air between them.

Inspired by the dynamics of computation, the authors of~\citep{vagner2014algebras} define a framework of open dynamical systems which compose as machines. When two machines compose, one machine is the designated sender and the other is the designated receiver. The sender emits information which directs the evolution of the receiver. In the special case where a system is both the sender and the receiver, this interaction describes feedback. When machines compose:
\begin{enumerate}
    \item \emph{communication is uni-directional}. Information travels from the sender to the receiver but not vice versa.
    \item \emph{interaction is active}. The communication channel from sender to receiver is specifically engineered to enable the passing of information. The receiver does not evolve without input from the sender. 
\end{enumerate}
When people communicate by passing notes, they are composing as machines where the note plays the role of the engineered communication channel.

The main theorem  of this paper (Theorem~\ref{thm:rsm}) unites these two flavors of composition in a single framework for open dynamical systems. In Section~\ref{sec:rsm-examples}, we define two concrete instances of dynamical systems --- continuous dynamical systems and non-deterministic automata --- and give examples of each composing as machines and as resource sharers. In Section~\ref{sec:rsm-operad}, we exemplify how operad and operad algebras respectively give a syntax and semantics for composition. This formalism will be the main tool we use to define compositions of open dynamical systems. In Section~\ref{sec:rsm-previous}, we discuss the established frameworks for composition as resource sharers and as machines in more depth. Finally, in Section~\ref{sec:rsm-main} we prove our main theorem.\sidenote{This paper contains sidenotes. Sidenotes contain two types of information. (1) Mathematical details for the sake of completeness but which detract from the main points of the paper. (2) "Recall that..." information intended to help the reader recall information presented earlier.

This document was formatted using the \texttt{tufte-latex} package in \LaTeX which I learned about from Tai-Danae Bradley.}

{\bf Acknowledgements} The author would like to thank David Spivak and David Jaz Myers for many helpful conversations.

\section{Motivating examples}\label{sec:rsm-examples}

The goal of this section is to present two examples of dynamical systems where both composition as machines and composition as resource sharers are fruitful methods of gluing dynamical systems together. 

\marginnote{{\bf Notation} For two morphisms $f: C \to D, g : D \to E$, we denote their composition $f \cp g: C \to E$.}

\subsection{Continuous Dynamical Systems}

A continuous dynamical system is defined by a state space $X$ and a vector field $v: X \to TX$.\sidenote{
Continuous dynamical systems often refers to a more general class of systems which are defined by a state space $X$ and a continuous group action $\Rb \curvearrowright X$ called a flow. However, in this paper we use the term "continuous dynamical systems" to refer to the subclass of systems with flows induced by a vector field.} In general, the state space $X$ may be any manifold, but in all of our examples $X$ will be a Euclidean space.
The vector field $v$ assigns to each $p\in X$ an arrow $v(p)$ based at $p$.\sidenote{
More formally, a vector field $v: X \to TX$ is a section of the tangent bundle $\pi: TX \to X$. So for $p \in X$, $v(p)$ is a vector in the tangent space $T_pX$.}.

The data of $(X,v)$ models "things that change" as follows: if the system is at a state $p \in X$, then it will evolve in the direction of $v(p)$. This intuition is particularly poignant if the state space $X$ represents the position of a ball. Then if the ball is at position $p \in X$, it will roll in the direction of the arrow $v(p)$.\sidenote{
For each $p \in X$ there is a unique trajectory $\gamma: \Rb \to X$ with velocity $\gamma'(t) = v(\gamma(t))$ and initial condition $\gamma(0) = p$. If a ball is dropped at position $p$, then after $t$ time  it will be at the position $\gamma(t)$.}

For example, the  growth of a population of rabbits can be modeled by the vector field $u: \Rb \to T\Rb$ defined by $$u(r) = \beta r \in T_r\Rb$$ where the state space $\Rb$ represents a population of rabbits.  If there are $r$ rabbits, then the rabbit population grows at a rate of $\beta r$.\sidenote{Given a starting rabbit population $r \geq 0$, the rabbit population will grow according to the unique trajectory $\gamma: \Rb \to \Rb$ with velocity $\gamma'(t) = u(\gamma(t))$ and initial condition $\gamma(0) = r$. So after $t$ time, the size of the rabbit population will be $\gamma(t)$.} This vector field may also be denoted by $\dot{r} = \beta r$.

Since we are interested in \emph{open} continuous dynamical systems, we will generalize to parameterized vector fields $v: X \times I \to TX$.\sidenote{Formally, a parameterized vector field $v: X \times I \to TX$ is a continuous map such that $v \cp \pi$ is projection on to the first coordinate, where $\pi: TX \to X$ is the natural projection map.} In this more general case, if a ball is in position $p \in X$ and it receives input $i \in I$, then it will roll in the direction of $v(p,i)$.

The growth of a fox population is parameterized by a population of prey that the foxes eat. This system is modeled by the parameterized vector field $v: \Rb \times \Rb \to T\Rb$ defined by $$v(f,e) = \alpha e f \in T_f\Rb$$ where  $f$ represents the fox population and  $e$ represents the population of prey to be eaten. This continuous dynamical system is equivalently denoted $\dot{f} = \alpha e f$. In this case, the fox population grows at a rate $\alpha$ according to the law of mass action.

Now, we will introduce two methods of composing continuous dynamical systems --- as machines and  as resource sharers.

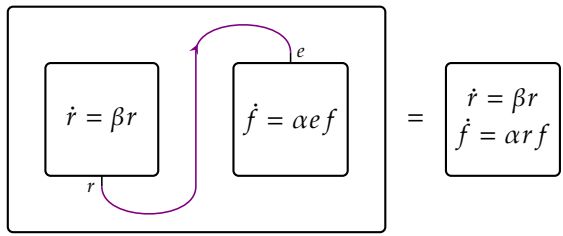
\begin{figure}\label{fig:rsm-cds-growth}
    \centering
    \begin{tikzpicture}[every text node part/.style={align=center}, oriented WD, bbx =1cm, bby =.5cm, bb min width=1.5cm, bb min height=1.5cm,bb port length=4pt, bb port sep=1]
	\node[bb={0}{0}{0}{1}] (X) {$\dot r = \beta r$};
	\node[bb={0}{0}{1}{0}, right = 1cm of X] (Y) {$\dot f = \alpha e f$};
	\node[bb={0}{0}{0}{0}, fit={($(X.north west)+(.5,.5)$) ($(Y.south east)+(-.5,-.5)$)}] (tot) {};
	
	\node[right = 1ex of tot](eq) {$=$};
	\node[bb={0}{0}{0}{0}, right= 1ex of eq] (Z) {$\dot r = \beta r$\\ $\dot f = \alpha rf$};
	\draw[label] 
		node [right = 2pt of Y_top1] {$e$}
		node [left = 2pt of X_bot1, yshift = -1ex] {$r$}
		;
		
     \draw[ar, color =violet] let \p1=(X.south east), \p2=(Y.north west), \n1=\bbportlen, \n2=\bby in
  		(X_bot1') to [out=-90, in=-90](\x1/2+\x2/2,\y1-\n1) -- (\x1/2+\x2/2,\y2+\n1) to [out=90,in=90](Y_top1);
\end{tikzpicture}
    \caption{This is an example of composing continuous dynamical systems as machines. Examining the left-hand side of the equation, the box on the left (the sender) is filled with a dynamical system modeling how a rabbit population grows. The box on the right (the receiver) is filled with a dynamical system modeling how a fox population grows  parameterized by an input population $e$. The directed wire indicates sending the rabbit population as  input  to the receiver. The total system (depicted on the right-hand side of the equation) models how the fox and rabbit populations grow in synchrony.}
\end{figure}

\newthought{First, how do  continuous dynamical systems} compose as \emph{machines}? Recall, the two dynamical systems we have introduced: $$\dot r = \beta r$$ modeling how a rabbit population grows and $$\dot f = \alpha e f$$ modeling how a fox population grows parameterized by a population $e$ of prey for the foxes to eat. Since foxes eat rabbits, we can send the rabbit population $r$  to the model for fox growth as the parameter $e$. Doing so composes the two systems. The resulting total system is then given by the vector field $\Rb^2 \to T\Rb^2$ defined by $$\dot{r} = \beta r, \quad \dot{f} = \alpha rf.$$ Figure~\ref{fig:rsm-cds-growth} depicts the composition of these systems as machines. 

As a second example (which will be used shortly), we can likewise compose as machines (1) the continuous dynamical system modeling how the fox population declines at a rate $\delta$ --- given by $\dot f = -\delta f$ and  (2) the continuous dynamical system modeling how the rabbit population declines as a rate $\gamma$ parameterized by a population $h$ of predators that hunt rabbits --- given by $\dot r = -\gamma hr$. The resulting total system $$\dot f = -\delta f,\quad \dot r = -\gamma fr$$ describes how both populations decline in synchrony.

\newthought{Second, how do continuous dynamical systems} compose as \emph{resource sharers}? Consider the two dynamical systems we have constructed: 
$$\dot{r} = \beta r,\quad \dot{f} = \alpha rf$$ modeling how the rabbit and fox populations grow and $$\dot f = -\delta f, \quad\dot r = -\gamma fr$$ modeling how they decline.

However there are not two separate rabbit populations, one that grows and one that declines. Rather both systems are referring to a shared pool of resources, in this case a population of rabbits. Likewise both systems are a referring to a shared population of foxes. 

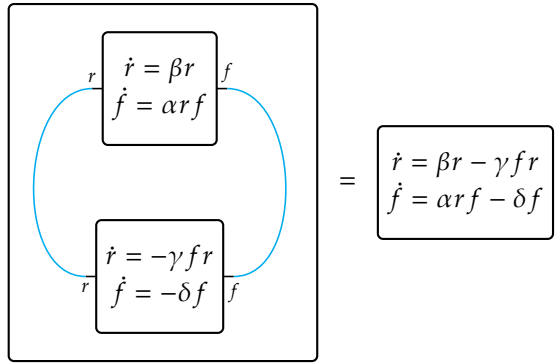
\begin{figure}\label{fig:CDS-rs}
    \centering
    \begin{tikzpicture}[every text node part/.style={align=center}, oriented WD, bbx = 1cm, bby =.5cm, bb min width=1.5cm, bb min height=1.5cm,bb port length=4pt, bb port sep=1]
	\node[bb={1}{1}{0}{0}] (B) {$\dot r = \beta r$\\ $\dot f = \alpha rf$};
	\node[bb={1}{1}{0}{0}, below=1cm of B] (D) {$\dot r = -\gamma f r$\\ $\dot f = -\delta f$};
	
	\draw[label] 
		node [above = 2pt of B_in1] {$r$}
		node [below = 2pt of D_in1] {$r$}
		node [above = 2pt of B_out1] {$f$}
		node [below = 2pt of D_out1] {$f$}
		;
    \draw[color=cyan] (B_in1) to [out=180, in=180] (D_in1);
    \draw[color=cyan] (B_out1) to [out=0, in=0] (D_out1);
    
    \node[bb={0}{0}{0}{0}, fit={($(B.north west)+(-.25,-.25)$) ($(D.south east)+(.25,.25)$)}] (tot) {};
    \node[right= 1ex of tot] (eq) {$=$};
    \node[bb={0}{0}{0}{0}, right=1ex of eq] (LV) {$\dot r = \beta r - \gamma fr$\\ $\dot f = \alpha rf - \delta f$};
\end{tikzpicture}
    \caption{This is an example of composing continuous dynamical systems as resource sharers. Examining the left-hand side of the equation, the top box is filled with a dynamical system modeling how the rabbit and fox populations grow. The bottom box is filled with a dyanmical system modeling how the rabbit and fox populations decline. The undirected wires connecting these boxes indicate that these two systems are composed by identifying their rabbit and fox populations respectively. The resulting dynamical system drawn on the right, is the Lokta-Volterra predator-prey model.}
\end{figure}

To compose these systems along the shared pools of rabbits and foxes, we add the effects of both systems on the shared resource. The resulting dynamical system is $$\dot{r} = \beta r - \gamma fr,\quad \dot {f} = \alpha rf - \delta f$$ known as the Lokta-Volterra predatory-prey model. Figure~\ref{fig:CDS-rs} depicts the composition of these systems as resource sharers.

\subsection{Non-deterministic Automata}

A non-deterministic automaton is a discrete dynamical system with a set of states $S$ and an update map $u: S \to \Pb(S)$.\sidenote[][4cm]{$\Pb$ denotes that power set monad, so $\Pb(S)$ is the set of subsets of $S$.}  For each state $s \in S$, the set of next possible states is $u(s) \subseteq S$.

\begin{figure}\label{fig:rsm-automata}
    \centering

    \subfloat[]{\begin{tikzpicture}[shorten >=1pt,node distance=2cm,on grid,auto] 
   \node[state] (q_0)   {$0$}; 
   \node[state] (q_1) [ right=of q_0] {$1$}; 
    \path[->] 
    (q_0) edge [bend left]  (q_1)
    (q_1) edge [bend left] node {\vphantom{1}} (q_0);
\end{tikzpicture}}\qquad
    \subfloat[]{\begin{tikzpicture}[shorten >=1pt,node distance=2cm,on grid,auto] 
   \node[state] (q_0)   {$0$}; 
   \node[state] (q_1) [ right=of q_0] {$1$}; 
    \path[->] 
    (q_0) edge [bend left] node{$1$}  (q_1)
        edge  [loop above] node {$0$} (q0)
    (q_1) edge [bend left] node{$1$} (q_0)
        edge  [loop above] node {$0$} (q1);
\end{tikzpicture}}

    \caption{(a) A 2-cycle automaton that oscillates between two states. \newline (b) A mod 2 adder automaton that takes as input $0$s and $1$s and adds the number of $1$s received modulo 2.}
\end{figure}
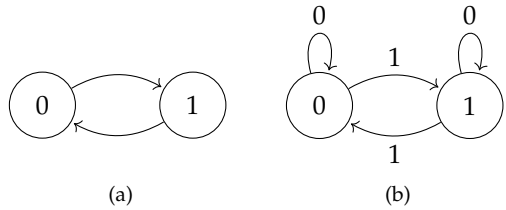

For example, consider an automaton representing a 2-cycle, depicted in Figure~\ref{fig:rsm-automata}(a). This automaton has two states and oscillates between them. Formally, $S = \Zb/2\Zb$ and $$u(s) = \{s + 1 \mod 2\}.$$

We want to generalize to \emph{open} automata, in other words automata that are parameterized by some input. These consist of a set of states $S$, a set of inputs $I$, and an update map $u: S \times I \to \Pb(S)$. For each state $s \in S$ and input $i \in I$, the set of next possible states is $u(s,i) \subseteq S$.

For example, consider an automaton that takes as input $0$s and $1$s and adds the number of $1$s received modulo 2, depicted in Figure~\ref{fig:rsm-automata}(b). Then $S = \Zb/2\Zb$, $I = \Zb/2\Zb$, and the update map is $$v(s,i) = \{s + i \mod 2\}.$$ 

Now, we will introduce two methods of composing non-deterministic automata --- as machines and  as resource sharers.

\newthought{First, how do automata} compose as  \emph{machines}? By reading off the states, the 2-cycle produces a string of $0$s and $1$s: $$... 01010101 ...$$ We can compose the automata from Figure~\ref{fig:rsm-automata} as machines by sending the string generated by the 2-cycle as input to the mod 2 adder. The output of the 2-cycle drives the dynamics of the mod 2 adder to define a total system.\sidenote[][-1in]{In Krohn-Rhodes theory, this style of composition is known as the cascade product of automata. See~\citep{krohn1965algebraic}.}

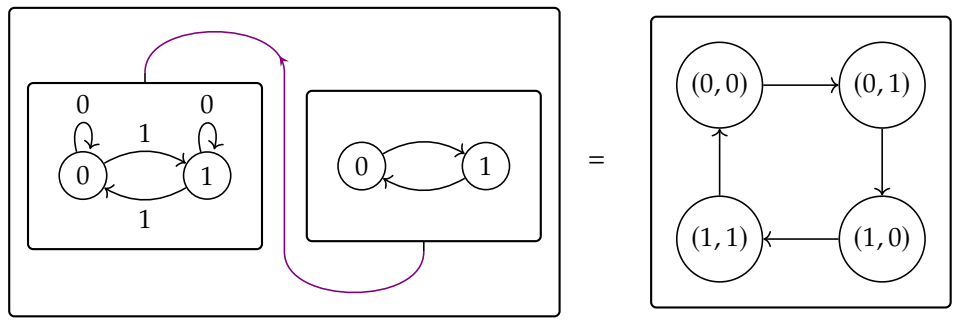
\begin{figure*}\label{fig:rsm-automata-machine}
    \centering
    \begin{tikzpicture}[every text node part/.style={align=center}, oriented WD, bbx = 1cm, bby =.5cm, bb min width=2cm, bb min height=2cm,bb port length=4pt, bb port sep=1]
   \node[state, minimum size=.5cm, minimum size=.5cm] (q_0)   {$0$}; 
   \node[state, minimum size=.5cm, minimum size=.5cm] (q_1) [ right=of q_0] {$1$}; 
    \path[->] 
    (q_0) edge [bend left] node[above]{$1$}  (q_1)
        edge  [loop above] node {$0$} (q_0)
    (q_1) edge [bend left] node[below] {$1$} (q_0)
        edge  [loop above] node {$0$} (q_1);
    
    \node[bb={0}{0}{1}{0}, fit={($(q_0.north west) + (.5, 1)$) ($(q_1.south east) + (-.5, -.5)$)}] (R) {};

    \node[state, minimum size=.5cm, right=1cm of R] (q2)   {$0$}; 
   \node[state, minimum size=.5cm] (q3) [ right=of q2] {$1$}; 
    \path[->] 
    (q2) edge [bend left]  (q3)
    (q3) edge [bend left] (q2);
    
    \node[bb={0}{0}{0}{1}, fit={($(q2.north west) + (.5, -.5)$) ($(q3.south east) + (-.5, .5)$)}] (S) {};
    
    \draw[ar, color=violet] let \p1=(S.south east), \p2=(R.north west), \n1=\bbportlen, \n2=\bby in
  		(S_bot1') to [out=-90, in=-90](\x1/2+\x2/2,\y1-\n1) -- (\x1/2+\x2/2,\y2+\n1) to [out=90,in=90](R_top1);
  		
  	\node[bb={0}{0}{0}{0}, fit = {($(R.north west)+(.75, 1)$) ($(S.south east) + (-.75, -1)$)}] (tot) {};

  	\node[right=2cm of  tot](x){};
  	\node[right= .25cm of  tot](){$=$};
  	
  	\node[state, minimum size=.5cm, above= 2ex of x](s1){$(0,0)$};
  	\node[state, minimum size=.5cm, below= 2ex of x](s4){$(1,1)$};
  	\node[state, minimum size=.5cm, right=of s1](s2){$(0,1)$};)
  	\node[state, minimum size=.5cm, right=of s4](s3){$(1,0)$};

  	\path[->]
  	(s1) edge (s2)
  	(s2) edge (s3)
  	(s3) edge (s4)
  	(s4) edge (s1);
  	
  	\node[bb={0}{0}{0}{0}, fit={($(s1.north west)+(.5, 0)$)($(s3.south east)+(-.5,0)$)}] (){};

\end{tikzpicture}
    \caption{ Composing the automata from Figure~\ref{fig:rsm-automata} as machines (left-hand side of the equation) yields a 4-cycle (right-hand side of the equation).}
\end{figure*}

The states of the total system are pairs of states of the individual systems so $\Zb/2\Zb \times \Zb/2\Zb$. Suppose we are in a state $(r,s)$ of the total system where $r$ belongs to the receiver and $s$ belongs to the sender. We update the total state by the following algorithm:
\begin{enumerate}
    \item send $s$ (the current state of the sender) to the receiver
    \item use the input to update the state of the receiver
    \item update the state of the sender
\end{enumerate} Formally, the update map is given by $$(r, s) \mapsto v(r,s) \times u(s) = \{(r + s \mod 2, s + 1 \mod 2)\}$$ and the total system is a 4-cycle. Figure~\ref{fig:rsm-automata-machine} depicts the composition of these automata as machines.

\newthought{Second, how do automata} compose as \emph{resource sharers}? In the case of automata, "resource sharing" represents "observation sharing" because it aligns automata along a shared observation. 

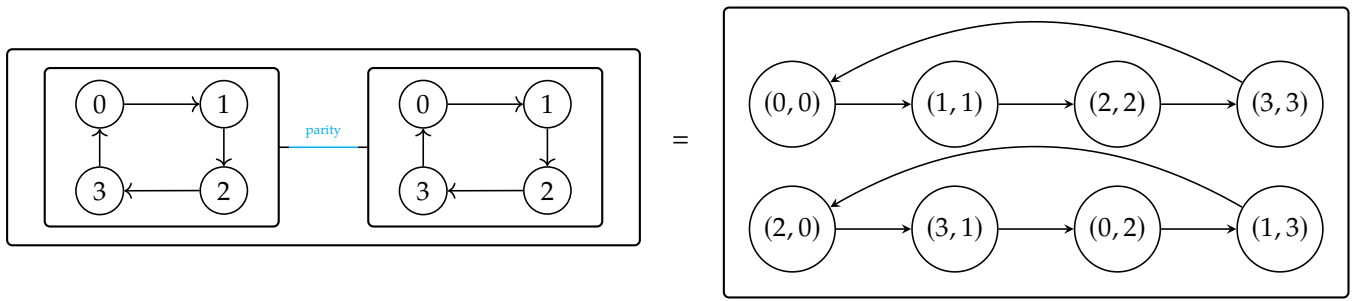
\begin{figure*}\label{fig:rsm-automata-rs}
    \centering
    \begin{tikzpicture}[every text node part/.style={align=center}, oriented WD, bbx = 1cm, bby =.5cm, bb min width=2cm, bb min height=2cm,bb port length=4pt, bb port sep=1]

  	\node[state, minimum size=.5cm](s1){$0$};
  	\node[state, minimum size=.5cm, right=of s1](s2){$1$};)
  	\node[state, minimum size=.5cm, below=of s2](s3){$2$};
  	\node[state, minimum size=.5cm, below=of s1](s4){$3$};
  	\path[->]
  	(s1) edge (s2)
  	(s2) edge (s3)
  	(s3) edge (s4)
  	(s4) edge (s1);
  	
  	\node[bb={0}{1}{0}{0}, fit={($(s1.north west) + (.5, -.5)$)($(s3.south east) -(.5, -.5)$)}] (c1){};
  	
  	\node[state, minimum size=.5cm, right=2cm of s2](t1){$0$};
  	\node[state, minimum size=.5cm, right=of t1](t2){$1$};)
  	\node[state, minimum size=.5cm, below=of t2](t3){$2$};
  	\node[state, minimum size=.5cm, below=of t1](t4){$3$};
  	\path[->]
  	(t1) edge (t2)
  	(t2) edge (t3)
  	(t3) edge (t4)
  	(t4) edge (t1);
  	
  	\node[bb={1}{0}{0}{0}, fit={($(t1.north west) +(.5, -.5)$)($(t3.south east)-(.5, -.5)$)}] (c2){};
  	
  	\draw[color=cyan] (c1_out1) edge node [above] {\tiny{parity}}  (c2_in1);

  	\node[bb={0}{0}{0}{0}, fit={($(c1.north west)+(.5, -.5)$)( $(c2.south east)-(.5, -.5)$)}] (tot){};
  	\node[right = 2.85cm of t1, yshift = -.5cm](eq){$=$};
  	
  	\node[state, minimum size=.5cm, right = 4cm of t1] (a1) {$(0,0)$};
  	\node[state, minimum size=.5cm, right = of a1] (a2) {$(1,1)$};
  	\node[state, minimum size=.5cm, right = of a2] (a3) {$(2,2)$};
  	\node[state, minimum size=.5cm, right = of a3] (a4) {$(3,3)$};
  	\path[-stealth]
  	(a1) edge (a2)
  	(a2) edge (a3)
  	(a3) edge (a4)
  	(a4) edge [bend right](a1);
  	
  	\node[state, minimum size=.5cm, below = of a1] (b1) {$(2,0)$};
  	\node[state, minimum size=.5cm, right = of b1] (b2) {$(3,1)$};
  	\node[state, minimum size=.5cm, right = of b2] (b3) {$(0,2)$};
  	\node[state, minimum size=.5cm, right = of b3] (b4) {$(1,3)$};
  	\path[-stealth]
  	(b1) edge (b2)
  	(b2) edge (b3)
  	(b3) edge (b4)
  	(b4) edge [bend right](b1);

  	\node[bb={0}{0}{0}{0}, fit={($(a1.north west)+(.5, .75)$)($( b4.south east) -(.5, 0)$)}] (){};

\end{tikzpicture}
    \caption{Two automata which count modulo $4$ sharing the observation of the parity of their states.}
\end{figure*}

Consider a 4-cycle which emits the parity of its states when observed.  We can align two such 4-cycles along the observation of state parity. The states of the total system are pairs of states of the individual 4-cycles which agree on parity. A transition is a pair of transitions of the individual 4-cycles where the domains agree on parity as do the codomains. In this example, the total system consists of two cycles representing the two ways  for the two 4-cycles to align along parity: either the states can agree exactly or they can be phase shifted by $2$. Figure~\ref{fig:rsm-automata-rs} depicts the composition of these automata as resource sharers.

We restrict our attention to non-deterministic automata because in some instances it is possible for a state in the total system to have no outgoing transitions.

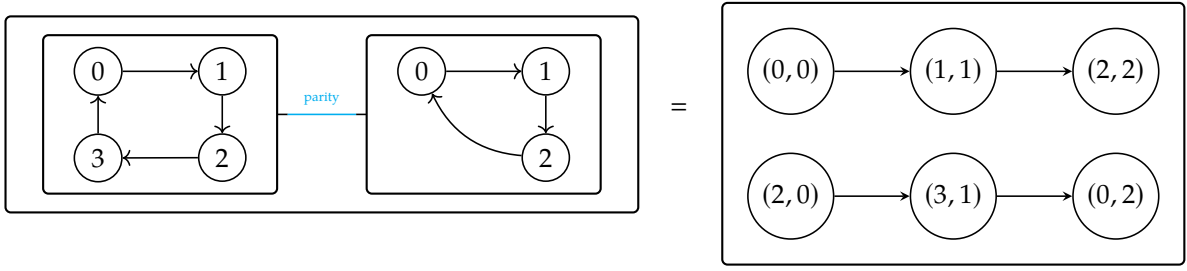
\begin{figure*}\label{fig:rsm-automata-rs2}
    \centering
    \begin{tikzpicture}[every text node part/.style={align=center}, oriented WD, bbx = 1cm, bby =.5cm, bb min width=2cm, bb min height=2cm,bb port length=4pt, bb port sep=1]

  	\node[state, minimum size=.5cm](s1){$0$};
  	\node[state, minimum size=.5cm, right=of s1](s2){$1$};)
  	\node[state, minimum size=.5cm, below=of s2](s3){$2$};
  	\node[state, minimum size=.5cm, below=of s1](s4){$3$};
  	\path[->]
  	(s1) edge (s2)
  	(s2) edge (s3)
  	(s3) edge (s4)
  	(s4) edge (s1);
  	
  	\node[bb={0}{1}{0}{0}, fit={($(s1.north west) + (.5, -.5)$)($(s3.south east) -(.5, -.5)$)}] (c1){};
  	
  	\node[state, minimum size=.5cm, right=2cm of s2](t1){$0$};
  	\node[state, minimum size=.5cm, right=of t1](t2){$1$};)
  	\node[state, minimum size=.5cm, below=of t2](t3){$2$};
  	\path[->]
  	(t1) edge (t2)
  	(t2) edge (t3)
  	(t3) edge [bend left] (t1);

  	\node[bb={1}{0}{0}{0}, fit={($(t1.north west) +(.5, -.5)$)($(t3.south east)-(.5, -.5)$)}] (c2){};
  	
  	\draw[color=cyan] (c1_out1) edge node [above] {\tiny{parity}}  (c2_in1);

  	\node[bb={0}{0}{0}{0}, fit={($(c1.north west)+(.5, -.5)$)( $(c2.south east)-(.5, -.5)$)}] (tot){};
  	\node[right = 2.85cm of t1, yshift = -.5cm](eq){$=$};
  	
  	\node[state, minimum size=.5cm, right = 4cm of t1] (a1) {$(0,0)$};
  	\node[state, minimum size=.5cm, right = of a1] (a2) {$(1,1)$};
  	\node[state, minimum size=.5cm, right = of a2] (a3) {$(2,2)$};
  	\path[-stealth]
  	(a1) edge (a2)
  	(a2) edge (a3);;
  	
  	\node[state, minimum size=.5cm, below = of a1] (b1) {$(2,0)$};
  	\node[state, minimum size=.5cm, right = of b1] (b2) {$(3,1)$};
  	\node[state, minimum size=.5cm, right = of b2] (b3) {$(0,2)$};
  	\path[-stealth]
  	(b1) edge (b2)
  	(b2) edge (b3);
  	
  	\node[bb={0}{0}{0}{0}, fit={($(a1.north west)+(.5, 0)$)($( b3.south east) -(.5, 0)$)}] (){};

\end{tikzpicture}
    \caption{A 4-cycle and a 3-cycle sharing the observation of state parity. Although the systems independently are deterministic automata, the total system is non-deterministic.}
\end{figure*}

For example, suppose we align along parity a 4-cycle and a 3-cycle. See Figure~\ref{fig:rsm-automata-rs2}. Consider the state $(2,2)$ of the total system. The $4$-cycle will update to its state $3$ while the $3$-cycle will update to its state $0$. Since these states differ in parity, there is no transition in the total system out of the state $(2,2)$.

\newthought{The goal of this paper} is to give an operad and operad algebra which describe composing dynamical systems as machines and as resource sharers \emph{simultaneously}. We will see the examples of continuous dynamical systems and non-deterministic automata as special cases of this formalism.

\section{Operad and Operad Algebras}\label{sec:rsm-operad}

In Section~\ref{sec:rsm-examples}, we graphically represented the composition of continuous dynamical systems and non-deterministic automata. Let's highlight the patterns in these pictures: (1) dynamical systems filled boxes and (2) composition corresponded to wiring boxes together. This section casually introduces the math behind those pictures, namely operads\sidenote{
In fact, this section is so casual that we only recite the definition of operad  for completeness!
\begin{defn}
An \emph{operad} $\Oa$ consists of a set of types $\ob \Oa$ and for types $s_1,.., s_n, t \in \ob \Oa$ a set of morphisms $\Oa(s_1,..., s_n; t)$ along with
\begin{itemize}
    \item for each type $t$, an identity morphism $\id_t \in \Oa(t; t)$
    \item a substitution map $\circ: \Oa(s_1,..., s_n; t_i) \times \Oa(t_1,..., t_i,..., t_m; u) \to \Oa(t_1,..., t_{i - 1}, s_1,..., s_n, t_{i + 1}, ..., t_m; u)$
    \item a symmetry map for each permutation of the domain types
\end{itemize}
satisfying an identity and associativity law. 
\end{defn} We refer the reader to~\citep{fong2019invitation} Chapter 6 for a helpful exposition of operads and to~\citep{leinster2004higher} for a complete definition. Note that this definition of an operad historically went by the name \emph{colored operad}.
} and operad algebras.

An operad $\Oa$ is much like a category. It has
\begin{itemize}
    \item a set of types $\ob \Oa$, analogous to objects in a category. Figure~\ref{fig:operad}(a) shows how we might visualize types as boxes.
    \item for types $s_1,...,s_n, t$, a set of morphisms\sidenote{Morphisms in an operads are sometimes refered to as \emph{operations}.}  $\Oa(s_1,...,s_n; t)$, analogous to morphisms in a category. Unlike in a category, an operad morphism  may have multiple (but finitely many) types as the domain. Figure~\ref{fig:operad}(b) shows how we might visualize morphisms in an operad as wirings between boxes.
\end{itemize} 

A symmetric monoidal category $(\Ca, \otimes, 1)$ induces an operad $\Oa(\Ca)$ with types $\ob \Ca$ and morphisms $$\Oa(\Ca)(s_1,..., s_n;t ) = \Ca(s_1 \otimes ... \otimes s_n, t).$$ All of the operads we will consider are induced by symmetric monoidal categories.

\begin{ex}
Consider the symmetric monoidal category $(\Finset\op, + , 0)$ and the induced operad $\Oa(\Finset\op)$. A type $M$ in $\Oa(\Finset\op)$ is a finite set. Graphically, a type $M$ is represented by a box with $M$ exposed ports. See Figure~\ref{fig:operad}(a). 

A morphism $f: (M_1, ..., M_n) \to N$ in $\Oa(\Finset\op)$ is a finite set map $f: N \to M_1 + ... + M_n$. Graphically, a morphism $f$ is represented by wiring each port $n \in N$ of the outer box to the port $f(n) \in M_i$ of an inner box. See Figure~\ref{fig:operad}(b). 
\end{ex}

\begin{figure}\label{fig:operad}
    \centering

    \subfloat[]{\begin{tikzpicture}[ oriented WD, bbx = 1cm, bby =.5cm, bb port length=4pt, bb port sep=.75]
    \node[bb={0}{6}{0}{0}] (){};
\end{tikzpicture}}\qquad
    \subfloat[]{\begin{tikzpicture}[ oriented WD, bbx = 1cm, bby =.5cm, bb port length=4pt, bb port sep=.75]
    \node[bb={0}{2}{0}{0}] (b1){};
    \node[bb={0}{3}{0}{0}, below= of b1] (b2){};
    \node[bb={0}{6}{0}{0}, fit={(b1.north west)(b2.south east)}] (tot){};
    
    \draw[color=cyan] (tot_out1') edge [out=180, in=0] (b1_out1);
    \draw[color=cyan] (tot_out2') edge [out=180, in=0] (b2_out1);
    \draw[color=cyan] (tot_out3') edge [out=180, in=0] (b1_out1);
    \draw[color=cyan] (tot_out4') edge [out=180, in=0] (b2_out1);
    \draw[color=cyan] (tot_out5') edge [out=180, in=0] (b2_out1);
    \draw[color=cyan] (tot_out6') edge [out=180, in=0] (b2_out3);

\end{tikzpicture}}

    \caption{(a) Graphically, we represent a type in an operad by a box with an interface. This box represents $6$ in the operad  $\Oa(\Finset\op)$. \newline
    (b) Graphically, we represent a morphism in $\Oa(s_1,..., s_n; t)$ as a wiring from the boxes corresponding to types $s_1,...,s_n$ to a box with type $t$. The inner boxes correspond to the domain types. The outer box corresponds to the codomain type. This wiring is a morphism $(2,3) \to 6$ in $\Oa(\Finset\op)$.}
\end{figure}

Operads give a syntax for composition, where a morphism in $\Oa(s_1,..., s_n; t)$ defines a way of composing types $s_1,..., s_n$ such that the result is of type $t$. Operads allow us to define a many different compositions between the same types. In analogy to the game MadLibs, types are like the parts of speech (noun, verb, adjective) and morphisms are ways of combining the parts of speech. There are many ways to compose two nouns and a verb into a sentence. For example, the MadLibs 
\begin{center}
    The $\underset{\mathsf{noun}}{\rule{1cm}{0.15mm}}$ in the $\underset{\mathsf{noun}}{\rule{1cm}{0.15mm}}$ $\underset{\mathsf{verb}}{\rule{1cm}{0.15mm}}$ed.
\end{center}
and
\begin{center}
    I $\underset{\mathsf{verb}}{\rule{1cm}{0.15mm}}$ed with joy when I got a $\underset{\mathsf{noun}}{\rule{1cm}{0.15mm}}$ and a $\underset{\mathsf{noun}}{\rule{1cm}{0.15mm}}$.\\
\end{center}
are two of many morphisms $(\mathsf{noun}, \mathsf{noun}, \mathsf{verb}) \to \mathsf{sentence}$ in the MadLibs operad.

To give meaning to a Madlibs, we must fill in each blank with a word of the appropriate type. The analogy continues. To give meaning to an operad, we must fill in each box with an element of the appropriate type. These semantics are given by the following structure.

\begin{defn}
Let $\Oa$ be an operad. An \emph{$\Oa$-algebra} is an operad functor\sidenote{ An operad functor $F: \Oa \to \Oa'$ consists of
\begin{itemize}
    \item a map of types $F: \ob \Oa \to \ob \Oa'$
    \item a map of morphisms $F: \Oa(s_1,..., s_n; t) \to \Oa'(Fs_1,..., Fs_n; Ft)$
\end{itemize} 
respecting identity and composition. Again we refer readers to~\citep{fong2019invitation} Chapter 6 for a helpful description and to~\citep{leinster2004higher} for a complete definition.
}
\begin{center}
    $F: \Oa \to \Oa(\Set, \times, 1)$.\sidenote{
    Recall that $\Oa(\Set, \times, 1)$ is the operad induced by the symmetric monoidal category $(\Set, \times, 1)$. In particular, its types are sets.}
\end{center}
\end{defn}

 An $\Oa$-algebra $F$ gives meaning to the syntax defined by $\Oa$ as follows:
\begin{itemize}
    \item on types --- for a type $t$ of $\Oa$, $F(t)$ is a set and $x \in F(t)$ is an element of type $t$. 
    
    \item on morphisms --- for a morphism $f \in \Oa(s_1,..., s_n; t)$ the set map $$Ff: F(s_1) \times ... \times F(s_n) \to F(t)$$ determines how composing elements of types $s_1,..., s_n$ according to $f$ results in an element of type $t$.
\end{itemize}

A lax monoidal functor $F: (\Ca, \otimes, 1) \to (\Set, \times, 1)$ induces an $\Oa(\Ca)$-algebra, which we occasionally refer to as a $\Ca$-algebra.

Suppose an operad $\Oa$ can be graphically represented by boxes and wirings.\sidenote{We have seen that $\Oa(\Finset\op)$ is an example of such an operad. Its graphical representation is shown in Figure~\ref{fig:operad}.} Then we represent an $\Oa$-algebra $F$ as follows. Elements of $F(t)$ fill boxes of type $t$. Let $f \in \Oa(s_1,..., s_n; t)$ be a wiring. Then the set map $Ff$ defines how filling the inner boxes with elements of the appropriate type and composing along the wiring defined by $f$ induces a filling for the outer box.

\begin{ex}
Given an alphabet $\Sigma$, there exists a lax monoidal functor $\mathsf{Label}_{\Sigma}: (\Finset\op, + , 0) \to (\Set, \times, 1)$ defined
\begin{itemize}
    \item on objects --- $M$ maps to $\Sigma^M$, the collection of labelings  $\sigma: M \to \Sigma$.
    \item on morphisms --- $f: M \to N$ maps to $f^*: \Sigma^N \to \Sigma^M$ defined by $f^*(\sigma)(m) = \sigma(f(m))$.
\end{itemize}

Under the induced $\Oa(\Finset\op)$-algebra, a box of type $M$ is filled with an element of $\Sigma^M$, i.e. a label in $\Sigma$ for each port $m \in M$. See Figure~\ref{fig:rsm-operad-labeled}(a). 

A finite set map $f: N \to M_1 + ... + M_n$ defines a wiring from  inner boxes of types $M_1,..., M_n$ to an outer box of type $N$. The set map $$\mathsf{Label}_{\Sigma}(f): \mathsf{Label}_{\Sigma}(M_1 + ... + M_n) \to \mathsf{Label}_{\Sigma}(N)$$ defines how filling each inner box with a choice in $\mathsf{Label}_{\Sigma}(M_i)$ and wiring along $f$ results in a labeling of type $N$. See Figure~\ref{fig:rsm-operad-labeled}(b).
\end{ex}
    
\begin{figure}\label{fig:rsm-operad-labeled}
    \centering

    \subfloat[]{\begin{tikzpicture}[ oriented WD, bbx = 1cm, bby =.5cm, bb port length=4pt, bb port sep=.75, bb min width=1.5cm, bb min height=1.5cm,]
    \node[bb={0}{3}{0}{0}] (b){};
    
    \node[left =3pt of  b_out1] {$b$};
    \node[left =3pt of  b_out2] {$a$};		
    \node[left =3pt of  b_out3] {$a$};		

\end{tikzpicture}}\qquad
    \subfloat[]{\begin{tikzpicture}[ oriented WD, bbx = 1cm, bby =.5cm, bb port length=4pt, bb port sep=.75]
    \node[bb={0}{2}{0}{0}] (b1){};
    \node[bb={0}{3}{0}{0}, below= of b1] (b2){};
    \node[bb={0}{6}{0}{0}, fit={(b1.north west)(b2.south east)}] (tot){};
    
    \draw[color=cyan] (tot_out1') edge [out=180, in=0] (b1_out1);
    \draw[color=cyan] (tot_out2') edge [out=180, in=0] (b2_out1);
    \draw[color=cyan] (tot_out3') edge [out=180, in=0] (b1_out1);
    \draw[color=cyan] (tot_out4') edge [out=180, in=0] (b2_out1);
    \draw[color=cyan] (tot_out5') edge [out=180, in=0] (b2_out1);
    \draw[color=cyan] (tot_out6') edge [out=180, in=0] (b2_out3);
    
    \node[left =3pt of  b2_out1] {$b$};
    \node[left =3pt of  b2_out2] {$a$};		
    \node[left =3pt of  b2_out3] {$a$};
    \node[left =3pt of  b1_out1] {$a$};
    \node[left =3pt of  b1_out2] {$b$};	
    
    \node[right= 2pt of tot] (eq) {$=$};
    \node[bb={0}{6}{0}{0}, right=2pt of eq, bb min height=3.5cm, bb min width= 1cm] (f){};

    \node[left=3pt of f_out1] {$a$};
    \node[left=3pt of f_out2] {$b$};
    \node[left=3pt of f_out3] {$a$};
    \node[left=3pt of f_out4] {$b$};
    \node[left=3pt of f_out5] {$b$};
    \node[left=3pt of f_out6] {$a$};

\end{tikzpicture}}

    \caption{(a) Suppose $\Sigma = \{a,b\}$. Graphically, a box of type $M$ is filled with a choice of labeling in $\Sigma^M$. Here we have a box of type $3$ that is filled with the labeling $(b,a,a) \in \Sigma^3$. \newline (b) The wiring depicted is an operad morphism $f: (2, 3) \to 6$. Filling the inner boxes with a labeling induces a labeling for the outer box. The label for each port of the outer box is determined by following the wire to a labeled inner port.}
\end{figure}

\newthought{In the remainder of this section} we will discuss  two  operads, which define syntaxes for composing dynamical systems as resource sharers and as machines.

\begin{ex}[Syntax for resource sharers]
Consider the symmetric monoidal category $(\Cospan_{\Finset}, + , 0)$. The induced operad $\Oa(\Cospan_{\Finset})$ has 
\begin{itemize}
    \item types --- finite sets $M$
    \item morphisms --- cospans $M_1 + ... + M_n \to Q \from N$
\end{itemize}

As in $\Oa(\Finset\op)$, a type $M$ is graphically represented by a box with $M$ exposed ports. To graphically represent a morphism $$M_1 + ... + M_n \xrightarrow{i} Q \xleftarrow{j} N$$ we draw an intermediate box with $Q$ exposed ports. Then we wire the ports from the inner boxes (respectively, outer box) to the intermediate box according to $i$ (respectively, $j$).  Often, we do not draw the intermediate box and simply draw the wiring of ports which may include (1) combining many ports irrespective of their origin and (2) terminating ports. See Figure~\ref{fig:rsm-cospan}.

\begin{figure}\label{fig:rsm-cospan}
    \centering

    \subfloat[]{\begin{tikzpicture}[ oriented WD, bbx = 1cm, bby =.5cm, bb port length=4pt, bb port sep=.75]
    \node[bb={0}{2}{0}{0}] (b1){};
    \node[bb={0}{3}{0}{0}, below= of b1] (b2){};
    \node[bb={0}{6}{0}{0}, fit={($(b1.north west)+(.75, 0)$)(b2.south east)}, color=lightgray] (q){};
    \node[bb={0}{5}{0}{0}, fit={($(q.north west)+(.75, 0)$)(q.south east)}] (m){};

    \draw[color=cyan] (b1_out1) edge [out=0, in=180] (q_out2');
    \draw[color=cyan] (b1_out2) edge [out=0, in=180] (q_out1');
    \draw[color=cyan] (b2_out1) edge [out=0, in=180] (q_out2');
    \draw[color=cyan] (b2_out2) edge [out=0, in=180] (q_out4');
    \draw[color=cyan] (b2_out3) edge [out=0, in=180] (q_out5');
    
    \draw[color=cyan] (m_out1') edge [out=180, in=0] (q_out2);
    \draw[color=cyan] (m_out2') edge [out=180, in=0] (q_out3);
    \draw[color=cyan] (m_out3') edge [out=180, in=0] (q_out3);
    \draw[color=cyan] (m_out4') edge [out=180, in=0] (q_out6);
    \draw[color=cyan] (m_out5') edge [out=180, in=0] (q_out5);

\end{tikzpicture}}\qquad
    \subfloat[]{\begin{tikzpicture}[ oriented WD, bbx = 1cm, bby =.5cm, bb port length=4pt, bb port sep=.75]
    \node[bb={0}{2}{0}{0}] (b1){};
    \node[bb={0}{3}{0}{0}, below= of b1] (b2){};
    \node[bb={0}{5}{0}{0}, fit={($(b1.north west)+(.75, 0)$)(b2.south east)}] (m){};

	\node[left=.35 of m_out1]  (link) {};
	\draw[color=cyan] (b1_out1) edge [out=0, in=180] (link.center);
	\draw[color=cyan] (b2_out1) edge [out=0, in=-90] (link.center);
	\draw[color=cyan] (m_out1') edge [out=180, in=0] (link.center);
	
	\draw[color=cyan] (m_out2') edge [out=180, in=180] (m_out3');
	
	\node[right=.1 of b2_out2]  (link2) {$\bullet$};
	\draw[color=cyan] (b2_out2) -- (link2.center);
	
	\node[right=.1 of b1_out2]  (link21) {$\bullet$};
	\draw[color=cyan] (b1_out2) -- (link21.center);
	
	\node[left=.1 of m_out4']  (link3) {$\bullet$};
	\draw[color=cyan] (m_out4') -- (link3.center);
	
	\draw[color=cyan] (b2_out3) edge[out=0, in=180] (m_out5');

\end{tikzpicture}}

    \caption{(a) This wiring is a morphism $2 + 3\to 6 \from 5$ in the operad $\Oa(\Cospan_{\Finset})$. 
    \newline 
    (b) This wiring is a simplified visualization of the wiring depicted in (a).}
\end{figure}

Since $\Cospan$-algebras are 1-equivalent to hypergraph categories~\citep{fong2018hypergraph} and  hypergraph categories are input-output agnostic,  $\Cospan_{\Finset}$ is a sensible syntactic setting for resource sharing.
\end{ex}

The operadic setting for machines requires a bit of setup, which is described in more detail in~\citep{schultz2016dynamical}.

\begin{defn}
Let $\Ca$ be a category. The category of $\Ca$-typed finite sets, $\mathsf{TFS}_{\Ca}$, has
\begin{itemize}
    \item objects --- pairs $(M \in \Finset, \tau: M \to \ob \Ca)$
    \item morphisms --- a morphism $f: (M, \tau) \to (M', \tau')$ is a map $f: M \to M'$ such that $f \cp \tau' = \tau$.
\end{itemize}

The category $\mathsf{TFS}_{\Ca}$ has a symmetric monoidal product given on objects by $$(M, \tau) +  (M', \tau') = (M + M', [\tau, \tau']: M + M' \to \ob \Ca)$$  with monoidal unit $(0, ! : 0 \to \ob \Ca)$.
\end{defn}

\begin{defn}
Let $\Ca$ be a category. There exists a symmetric monoidal category of $\Ca$-wiring diagrams, $\Wd_{\Ca}$,  with
\begin{itemize}
    \item objects --- pairs of $\Ca$-typed finite sets $\lens{X^{\mathsf{in}}}{X^{\mathsf{out}}}$
    \item morphisms\sidenote{The morphisms in this category are often called prisms.} --- pairs $\lens{\phi^{\mathsf{in}}}{\phi^{\mathsf{out}}}: \lens{X^{\mathsf{in}}}{X^{\mathsf{out}}} \rightleftarrows \lens{Y^{\mathsf{in}}}{Y^{\mathsf{out}}}$ where $$\phi^{\mathsf{out}}: Y^{\mathsf{out}} \to X^{\mathsf{out}}, \quad \phi^{\mathsf{in}}: X^{\mathsf{in}} \to X^{\mathsf{out}} + Y^{\mathsf{in}}$$ are morphisms of $\Ca$-typed finite sets.
\end{itemize}
The symmetric monoidal structure on $\Wd_{\Ca}$ is induced by the symmetric monoidal strucutre of $\mathsf{TFS}_{\Ca}$.
\end{defn}

\begin{ex}[Syntax for machines]
The operad induced by $\Wd_{\Ca}$ is an appropriate syntax for composing dynamical systems as machines.

For $\Ca$-typed finite sets 
\begin{center}
    $X^{\mathsf{in}} = (M^{\mathsf{in}}, \tau^{\mathsf{in}}), \quad X^{\mathsf{out}} = (M^{\mathsf{out}}, \tau^{\mathsf{out}})$\sidenote{Recall that $M^{\mathsf{in}}$ is a finite set and $\tau^{\mathsf{in}}: M^{\mathsf{in}} \to \ob \Ca$ assigns to each port $m$ of $M^{\mathsf{in}}$ an object of $\Ca$. Likewise for $M^{\mathsf{out}}$ and $\tau^{\mathsf{out}}$.}
\end{center} we interpret the type $\lens{X^{\mathsf{in}}}{X^{\mathsf{out}}}$ of $\Oa(\Wd_{\Ca})$ as having $M^{\mathsf{in}}$ input wires where the wire $m \in M^{\mathsf{in}}$ carries information of type $\tau^{\mathsf{in}}(m)$. Likewise for the output wires.  See Figure~\ref{fig:rsm-wd}(a). We interpret morphisms $$\lens{\phi^{\mathsf{in}}}{\phi^{\mathsf{out}}}: \lens{X_1^{\mathsf{in}}}{X_1^{\mathsf{out}}} + ... + \lens{X_n^{\mathsf{in}}}{X_n^{\mathsf{out}}} \rightleftarrows \lens{Y^{\mathsf{in}}}{Y^{\mathsf{out}}},$$ as wiring diagrams like the one shown in Figure~\ref{fig:rsm-wd}(b).

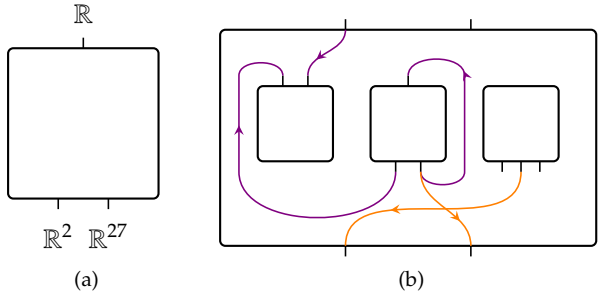
\begin{figure}\label{fig:rsm-wd}
    \centering

    \subfloat[]{\begin{tikzpicture}[oriented WD, bbx = 1cm, bby =.5cm, bb port length=4pt, bb port sep=.75, bb min width=2cm, bb min height=2cm]
    \node[bb={0}{0}{1}{2}] (B){};
    
    \node[above=2pt of B_top1] () {$\Rb$};
    \node[below=2pt of B_bot2'] () {$\Rb^{27}$};
    \node[below=2pt of B_bot1'] () {$\Rb^{2}$};

\end{tikzpicture}}\qquad
    \subfloat[]{\begin{tikzpicture}[oriented WD, bbx = .5cm, bby =.25cm, bb port length=4pt, bb port sep=.75, bb min width=1cm, bb min height=1cm]
    \node[bb={0}{0}{2}{0}] (A){};
    \node[bb={0}{0}{1}{2}, right=of A] (B){};
    \node[bb={0}{0}{0}{3}, right=of B] (C){};
    
    \node[bb={0}{0}{2}{2}, fit={($(A.north west)+(0,2)$)($(C.south east)+(0, -3.5)$)}] (tot){};
    
    \draw[ar, color=violet] (tot_top1') to [out=-90, in=90] (A_top2);
    
    \draw[ar, color=violet] let \p1=(A.south west), \p2=(A.north west), \n1=\bbportlen, \n2=\bby in
     (B_bot1') to [out=-90, in=-90] (\x1 - \bby, \y1 - \n1) to [out=90, in=-90] (\x2 - \bby, \y2 +\n1) to [out=90, in=90] (A_top1);
    
    \draw[ar, color=violet] let \p1=(B.south east), \p2=(B.north east), \n1=\bbportlen, \n2=\bby in
     (B_bot2') to [out=-90, in=-90] (\x1 + \bby, \y1 - \n1) to [out=90, in=-90] (\x2 + \bby, \y2 +\n1) to [out=90, in=90] (B_top1);
     
    \draw[ar, color=orange] (B_bot2') to [out=-90, in=90] (tot_bot2);
    \draw[ar, color=orange] (C_bot2') to [out=-90, in=90] (tot_bot1);
    
    \node[above=2pt of tot_top1'] () {};
    \node[above=2pt of tot_top2'] () {};
    
    \node[below=2pt of A_top1'] () {};
    \node[below=2pt of A_top2'] () {};
     \node[below=2pt of B_top1'] () {};

\end{tikzpicture}}

    \caption{(a) A type in $\Wd_{\Ca}$ is visualized as a box with a finite number of input and output ports labeled by objects of $\Ca$. This box represents a type in $\Wd_{\Euc}$ where $\Euc$ is the category of Euclidean spaces. \newline (b) A morphism in $\Wd_{\Ca}$ is visualized as two sets of wires. The purple wires, representing $\phi^{\mathsf{in}}$, feed the inputs to the inner boxes with either (1) outputs of the inner boxes or (2) inputs of the outer box. The orange wires, representing $\phi^{\mathsf{out}}$, feed the outputs of the outer box with outputs of the inner boxes. In this figure, the labeling of the ports has been suppressed however ports connected by a wire must be labeled with the same object of $\Ca$.}
\end{figure}

\end{ex}

\section{Previous Work}\label{sec:rsm-previous}

In this section we discuss two bodies of work, one which introduces resource sharers as a $\Cospan_{\Finset}$-algebra and one which introduces machines as a $\Wd_{\Ca}$-algebra. The motivation for the main theorem of this paper is uniting these two perspectives.

\newthought{In~\citep{Baez_2017}, the authors} define a hypergraph category $\Dynam$ which describes the composition of continuous dynamical systems as resource sharers. Recall that hypergraph categories are 1-equivalent to $\Cospan$-algebras. Taking the operadic perspective, $\Dynam$ is an $\Oa(\Cospan_{\Finset})$-algebra.  We consider each type $M$ to be a finite set of exposed ports.  The set $\Dynam(M)$ consists of triples $$(S \in \Finset, v: \Rb^S \to \Rb^S, p: S \to M)$$ where $v$ is an algebraic vector field\sidenote{
A vector field $v$ is algebraic if its components are polynomials.
} and $p$ is a map of finite sets. Therefore, under $\Dynam$, boxes with $M$ exposed ports are filled with triples $(S, v, p) \in \Dynam(M)$ where 
\begin{itemize}
    \item $S$ is a finite set of   state variables
    \item $v$ gives the dynamics on the state space $\Rb^S$
    \item $p$ assigns a state variable to each exposed port 
\end{itemize}

The details of how $\Dynam$ acts on morphisms in $\Cospan_{\Finset}$ are subsumed by the discussion following Theorem~\ref{thm:rsm}. Some intuition is given through the example shown in Figure~\ref{fig:CDS-rs}.

\newthought{In \citep{schultz2016dynamical}, the authors} define a $\Wd_{\Euc}$-algebra\sidenote{The category $\mathsf{Euc}$ is the full subcategory of $\Mnfld$ generated by Euclidean spaces.} $\CDS$\sidenote{
$\CDS$ stands for "continuous dynamical system."}  as follows.

Let $\Rb^I$ and $\Rb^O$ be Euclidean spaces. An $\lens{\Rb^I}{\Rb^O}$ continuous dynamical system is a triple $$(\Rb^S \in \Euc, v: \Rb^S \times \Rb^I \to T\Rb^S, r: \Rb^S \to \Rb^O)$$ where $v$ is a parameterized vector field.\sidenote[][-1cm]{
Recall that for a manifolds $X$ and $Y$, $v: X \times Y \to TX$ is a parameterized vector field on $X$ if the diagram below commutes. The map $TX \to X$ is the natural projection map.
}\begin{marginfigure}
\begin{center}
\begin{tikzcd}
X \times Y \arrow[r, "v"] \arrow[d, "\pi_1", swap]
    & TX \arrow[ld]\\
X & \\
\end{tikzcd}
\end{center}
\end{marginfigure}

There exists an algebra $\CDS: \Wd_{\mathsf{Euc}} \to \Set$ which on objects maps $\lens{X^{\mathsf{in}}}{X^{\mathsf{out}}}$\sidenote{
Recall that $\lens{X^{\mathsf{in}}}{X^{\mathsf{out}}} \in \ob \Wd_{\Euc}$ is a pair of $\Euc$-typed finite sets. We visualize $\lens{X^{\mathsf{in}}}{X^{\mathsf{out}}}$ as a box with a finite number of input and output ports each labeled with a Euclidean space.
} to the set of $\lens{\widehat{X^{\mathsf{in}}}}{\widehat{X^{\mathsf{out}}}}$ continuous dynamical systems where for a $\Euc$-typed finite set $X = (M, \tau: M \to \Euc)$, 
\begin{center}
    $\widehat{X} = \prod_{m \in M} \tau(m) \in \Euc$.\sidenote{
In general if $\Ca$ has finite products then there exists a functor $\widehat{(-)}: \mathsf{TFS}(\Ca) \to \Ca$ defined by  product.
}
\end{center}

Therefore, under $\CDS$, boxes of type $\lens{\Rb^I}{\Rb^O}$ are filled with triples $(\Rb^S, v, r)$ which we interpret as having
\begin{itemize}
    \item a state space --- the Euclidean space $\Rb^S$
    \item dynamics --- the vector field $v$ of $\Rb^S$ parameterized by the input space $\Rb^I$
    \item read-out --- the map $r$ taking states of $\Rb^S$ to points in the output space $\Rb^O$
\end{itemize}

The details of how $\CDS$ acts on morphisms are subsumed by the discussion following Theorem~\ref{thm:rsm}. Some intuition is given through the example shown in Figure~\ref{fig:rsm-cds-growth}.

\newthought{Notice the strong similarities} between the sets $\Dynam(M)$ and $\CDS\left(\lens{\Rb^I}{\Rb^O}\right)$. Both contain triples which determine (1) a state space, (2) dynamics, (3) an observation of the state space. However, the two algebras define remarkably different compositions of dynamical systems. In the next section, we give a single algebra capturing both types of composition.

\newpage

\section{Resource Sharing Machines}\label{sec:rsm-main}

To define resource sharing machines we will need two mathematical tools, which together define the data of a \emph{contravariant dynamical system doctrine}. The first tool is a category of lenses whose morphisms capture the data of an open dynamical system with a machine-style interface. The key observation: \begin{center}\emph{The data of an open dynamical system is given by a lens.}
\end{center} is made in \citep{spivak2019generalized} and further explored in \citep{myers2020double}.

\begin{defn}
For an indexed  category\sidenote{$\Aa: \Ca\op \to \Cat$ is a (monoidal) indexed category means that $\Aa$ is a (lax monoidal) psuedofunctor.}  $\Aa: \Ca \op \to \Cat$ define $$\Lens_{\Aa} := \int^{C \in \Ca} \Aa(\Ca)\op,$$ the pointwise opposite of the Grothendieck construction.
\end{defn}

Let's unpack this definition. The category $\Lens_{\Aa}$ has
\begin{itemize}
    \item objects --- pairs $\lens{I}{O}$ for $O \in \ob \Ca$ and $ I \in \ob \Aa(O)$
    \item morphisms --- lenses $$\lens{f^\#}{f}: \lens{I}{O} \leftrightarrows \lens{I'}{O'}$$ with $f: O \to O'$ in $\Ca$ and $f^\#: f^*I' \to I$ in $\Aa(O)$\sidenote{For $f: O \to O'$ in $\Ca$ we denote $\Aa(f)$ by $f^*$.}
\end{itemize}

The key observation is that for well-chosen indexed categories $\Aa$, certain lenses
\begin{equation*}
    \lens{u}{r}: \lens{TS}{S} \leftrightarrows \lens{I}{O}
\end{equation*}  capture the data of an open dynamical systems with
\begin{itemize}
    \item $S$ --- the internal state space
    \item $O$  ---  a space of outputs or orientations
    \item $I$  ---  a space of contextualized inputs over $O$
    \item $TS$  ---  a space of canonical changes over $S$
    \item $r: S \to O$ ---  a read-out map
    \item $u: r^*I \to TS$  ---  an update map
\end{itemize}

To see this observation in action, let's consider two examples. For each example, watch for (1) a choice of indexed category $\Aa$ and (2) a description of how specific lenses in $\Lens_{\Aa}$ correspond to open dynamical systems.

\begin{ex}[non-deterministic automata]
Let $\Pb: \Set \to \Set$ be the powerset monad. Consider the indexed  category

\begin{center}
$\biKl((-)  \times -, \Pb) : \Set \op \to \Cat$.\sidenote[][-.5in]{
For a set $S$, the category $\biKl(S  \times -, \Pb)$ has
\begin{itemize}
    \item objects --- sets, $X$
    \item morphisms --- $$\biKl(S  \times -, \Pb)(X, Y) = \Set(S \times X, \Pb Y)$$
\end{itemize}
On morphisms $f$, the functor $\biKl(f \times -, \Pb)$ is the identity on objects and  is  precomposition by $f \times \id$ on morphisms.
}
\end{center}

Unpacking definitions, a lens \begin{equation*}
    \lens{u}{r}: \lens{S}{S} \leftrightarrows \lens{I}{O}
\end{equation*} in $\Lens_{\biKl((-)  \times -, \Pb)}$ consists of
\begin{itemize}
    \item $S$ --- a set of states of the automaton
    \item $O$ --- a set of outputs of the automaton
    \item $I$ --- a set of inputs to the automaton
    \item $r: S \to O$ --- a set map assigning an output to each state
    \item $u: S \times I \to \Pb S$ --- a set map assigning to a state $s$ and input $i$, a set of next possible states\sidenote[][-.5in]{In this example, $TS$ (the space of canonical changes over $S$) is $S$ itself because a transition in an automaton is a choice of states in $S$ to transition to.}
\end{itemize}

\end{ex}

\begin{ex}[continuous dynamical systems]\label{ex:rsm-cds-generalized-lens}
Let $\Mnfld_{\Sub}$ be the wide subcategory of $\Mnfld$ whose morphisms are submersions.\sidenote{
We restrict our attention to submersions because in $\Mnfld$ pullbacks along submersions always exist.
}

Consider the indexed category $$\Mnfld_{\Sub}/(-): \Mnfld\op \to \Cat$$ which takes a manifold $B$ to the category of submersions over $B$.\sidenote{
$\Mnfld_{\Sub}/B$ is the category with
\begin{itemize}
    \item objects --- submersion $p: E \to B$
    \item morphisms --- commuting triangles
\end{itemize}
} \begin{marginfigure}
\begin{center}
\begin{tikzcd}
    E_1 \arrow[rr, "f"] \arrow[dr, "p_1", swap]
        && E_2 \arrow[dl, "p_2"]\\
        & B
        &\\
\end{tikzcd}
\end{center}
\end{marginfigure}
\marginnote{
For a continuous map $f: B \to B'$, the functor $\Mnfld_{\Sub}/f$ is given by taking pullbacks.}

Unpacking definitions, a lens
    $$\lens{u}{r}: \lens{TS \to S}{S} \leftrightarrows \lens{I \times O \xrightarrow{\pi_2} O}{O}$$
 in $\Lens_{\Mnfld_{\Sub}/(-)}$ consists of
\begin{itemize}
    \item $S$ --- a manifold giving the state space of the dynamical system
    \item $O$ --- a manifold of outputs
    \item $I$ --- a manifold of inputs
    \item $r: S \to O$ --- a continuous map assigning an output to each state of the state space
    \item $u: S \times I \to TS$ --- an indexed vector field assigning to a state $s$ and input $i$, a vector $u(s,i) \in T_sS$ indicating a direction in which to evolve
\end{itemize}

The dynamical systems described above are the special cases where the fiber of inputs over each output $o \in O$ is constantly $I$. In the general case, the inputs may vary with the output.\sidenote[][-1.25in]{ More generally, a lens
\begin{equation*}
    \lens{u}{r}: \lens{TS\to S}{S} \leftrightarrows \lens{p:I \to O}{O}
\end{equation*} in $\Lens_{\Mnfld_{\Sub}/(-)}$ consists of
\begin{itemize}
    \item $S$ --- a manifold giving the state space of the dynamical system
    \item $O$ --- a manifold of outputs
    \item $I$ --- a manifold of inputs where the fiber $p\inv(o) \subseteq I$ is the space of contextualized inputs for the output $o$
    \item $r: S \to O$ --- a continuous map assigning an output to each state of the state space
    \item $u: S \times_O I \to TS$ --- an indexed vector field assigning to each state $s$ and contextualized input $i \in p\inv(r(s))$, a vector $u(s,i) \in T_sS$ indicating a direction in which to evolve
\end{itemize}
}
\end{ex}

As shown in \citep{moeller2018monoidal}, if $\Aa: \Ca\op \to \Cat$ is a symmetric monoidal indexed category then $\Lens_{\Aa}$ has a symmetric monoidal structure,\sidenote{The monoidal structure for $\Lens_{\Aa}$ is given by  $$\lens{I}{O} \otimes \lens{I'}{O'} = \lens{I \boxtimes I'}{O \otimes O'}$$ where $O\otimes O'$ is the monoidal product in $\Ca$ and $I \boxtimes I'$ is the image of $(I, I') \in \Aa(O) \times \Aa(O')$ under the laxator.} and thus induces an operad $\Oa(\Lens_{\Aa})$.

\newthought{Morphisms of $\Lens_{\Aa}$ define dynamics} with a machine-style interface that distinguishes input and output. To define the resource sharing interface for resource sharing machines, we develop a second tool --- a section $T: \Ca \to \Lens_{\Aa}$ of the forgetful functor $U : \Lens_{\Aa} \to \Ca$.\sidenote[][-1in]{$U: \Lens_{\Aa} \to \Ca$ is defined
\begin{itemize}
    \item on objects --- $\lens{I}{O} \mapsto O$
    \item on morphisms --- $\lens{f^\#}{f} \mapsto f$
\end{itemize}
} Together an indexed category $\Aa: \Ca\op \to \Cat$ and a section $T: \Ca \to \Lens_{\Aa}$ is the data we need to specify a class of resource sharing machines. We refer to this data as a \emph{contravariant dynamical system doctrine}\sidenote{\emph{Contravariant dynamical system doctrines} stand in contrast to the (covariant) dynamical system doctrines defined in \citep{myers2020double} Definition 1.1.
\begin{itemize}
    \item A covariant dynamical system doctrine consists of an indexed category $\Aa: \Ca \op \to \Cat$ along with a section $T$ of its Grothendieck construction.
    \item In contrast, a contravariant dynamical system doctrine consists of an indexed category $\Aa: \Ca \op \to \Cat$ along with a section $T$ of the pointwise opposite of its Grothendieck construction.
\end{itemize} } or simply as a \emph{contravariant doctrine}.

\begin{defn}
A \emph{contravariant dynamical system doctrine} is an indexed category $\Aa: \Ca \op \to \Cat$ with a section $T: \Ca \to \Lens_{\Aa}$.
\end{defn}

Given a contravariant dynamical system doctrine $(\Aa: \Ca \op \to \Cat, T)$, we notate $T$
\begin{itemize}
    \item on objects, by $S \mapsto \lens{TS}{S}$
    \item on morphisms, by $f \mapsto \lens{Tf^\#}{f}$
\end{itemize}

Given a map of internal state spaces $f: S \to S'$, we say that $$Tf^\#: f^* (TS') \to TS$$ pullsback transitions over $S'$ to transitions over $S$ along $f$. To see this structure at play in the context of resource sharing, let's continue with our two examples.

\begin{ex}[non-deterministic automata]
Recall that the relevant indexed category for non-deterministic automata is $$\biKl((-) \times -, \Pb): \Set\op \to \Cat.$$ Define $T: \Set \to \Lens_{\biKl((-) \times -, \Pb)}$
\begin{itemize}
    \item on objects --- $S$ maps to $\lens{S}{S}$
    \item on morphisms --- $f: S \to S'$ maps to $\lens{Tf^\#}{f}: \lens{S}{S} \leftrightarrows \lens{S'}{S'}$ where the set map $Tf^\#: S \times S' \to \Pb S$ is defined by $$Tf^\#(s, s') = f\inv(s') \subseteq S.$$
\end{itemize}

How should we think about the lens $\lens{Tf^\#}{f}$? Interpret $f: S \to S'$ as a map between the states of two automata. If the second automaton transitions to state $s' \in S'$, then the first automaton must transition to a state in $f\inv(s')$.

The pair $(\biKl((-) \times -, \Pb), T)$ is a contravariant dynamical system doctrine.

\end{ex}

\begin{ex}[continuous dynamical systems] \label{ex:rsm-cds-section-T}
Recall that the relevant indexed category for continuous dynamical systems is $$\Mnfld_{\Sub}/(-): \Mnfld\op \to \Cat.$$ In order to define a contravariant dynamical system doctrine we must adjust this indexed category.

Let $\Riem$ be the category of Riemannian manifolds\sidenote{A Riemannanian manifold $(M, g)$ is a manifold $M$ equipped with a Riemannian structure $g$. The Riemannian structure allows us to define notions of length and angle on $M$. For our purposes $g$ defines an inner product on the tangent space $T_pM$ for each $p \in M$ that varies continuously with $M$. Importantly, $g$ induces a natural isomorphism between the tangent space $T_pM$ and the cotangent space $T_p^*M$. Let $\theta^M_p: T_pM \to T_p^*M$ denote this isomorphism.
}
with differentiable maps between them. In particular, a Euclidean space equipped with the standard inner product is an instance of a Riemannian manifold.

Composing $\Mnfld_{\Sub}/(-)$ with the functor $\Riem \to \Mnfld$ that forgets the Riemannian structure induces a new indexed category
\begin{equation*}
    \Riem_{\Sub}/(-): \Riem\op \to \Mnfld\op \xrightarrow{\Mnfld_{\Sub}/(-)} \Cat
\end{equation*}
An object of $\Lens_{\Riem_{\Sub}/(-)}$ is a pair $\lens{I \xto{p} O}{(O, g)}$ where $g$ is a Riemannian structure on $O$ and $p$ is a submersion of manifolds. Note that $I$ need not be Riemannian.

Define $T: \Riem \to \Lens_{\Riem_{\Sub}/(-)}$ as follows:
\begin{itemize}
    \item on objects --- $(S, g)$ maps to $\lens{TS \to S}{(S,g)}$ where $TS \to S$ is the natural projection map
    \item on morphisms --- $f: S \to S'$ maps to $$\lens{Tf^\#}{f}: \lens{TS \to S}{(S,g)} \leftrightarrows \lens{TS' \to S'}{(S', g')}$$ where $Tf^\#: TS' \times_{S'} S\to TS$ is defined  as follows. For $x \in S, \tilde y \in T_{f(x)}S'$, let $$T^\#f(\tilde y, x)  = ((\theta^S_x)\inv \circ T_x^*f \circ \theta_{f(x)}^{S'})(\tilde y).$$
\end{itemize}

This definition is slick but obscures the relationship to resource sharing. To achieve this more earthly goal, let's restrict our attention to Euclidean spaces.

Let $f: S' \to S$ be a map of finite sets. Then $f$ induces a map of Euclidean spaces $f^*: \Rb^{S} \to \Rb^{S'}$. The slogan for $T(f^*)^\#$ is \emph{add along shared coordinates}. For each $p \in \Rb^S$, the tangent space $T_p(\Rb^S)$ is isomorphic to  $\Rb^S$. We interpret elements of $T_p\Rb^S$ as maps $v: S \to \Rb$ which assign to each coordinate $s \in S$ a velocity $v(s)$. Under this interpretation, $T\Rb^{S'} \times_{\Rb^{S'}} \Rb^S \iso \Rb^{S'} \times \Rb^S$ and $$T(f^*)^\#(v, p) = \sum_{s' \in f\inv(p)} v(s')$$ justifying the slogan. If two coordinates in $S'$ are identified by $f$, then their velocities are summed in $T(f^*)^\#(v,p)$.

The pair $(\Riem_{\Sub}/(-), T)$ is a contravariant dynamical system doctrine.

\end{ex}

\newthought{We are now ready} to use the two tools provided by a contravariant dynamical system doctrine --- the indexed category $\Aa: \Ca\op \to \Cat$ and the section $T: \Ca \to \Lens_{\Aa}$ --- to define resource sharing machines.

\begin{thm}\label{thm:rsm}
Let $\Ca$ be a cartesian monoidal category with pullbacks. Let $(\Aa: \Ca\op \to \Cat, T: \Ca \to \Lens_{\Aa})$ be a contravariant dynamical system doctrine such that $\Aa$ is monoidal and $T$ is oplax. There exists a lax monoidal functor $$\Rsm: \Lens_{\Aa} \times \Span_{\Ca} \to \Set$$ defined
\begin{itemize}
    \item on objects --- for $\lens{I}{O} \in \ob \Lens_{\Aa}$ and $M \in \ob \Ca$ $$\Rsm\left(\lens{I}{O}, M\right) = \left\{\left( S, \lens{u}{r}, p \right) \big| S \in \ob \Ca, \lens{u}{r} \in \Lens_{\Aa}\left(\lens{TS}{S} , \lens{I}{O}\right), p \in \Ca (S, M) \right\}.$$

    \item on morphisms --- For $\lens{f^\#}{f}: \lens{I}{O} \leftrightarrows \lens{I'}{O'}$ in $\Lens_{\Aa}$ and span $M \leftarrow Q \rightarrow M'$ in $\Span_{\Ca}$, the set map
    $\Rsm\left(\lens{f^\#}{f}, M \xleftarrow{i} Q \xrightarrow{i'} M' \right)$ maps the triple $\left(S, \lens{u}{r}, p\right) \in \Rsm\left(\lens{I}{O}, M\right)$ to $$\left(S \times_M Q, \lens{T\tilde i ^\#}{\tilde i } \cp \lens{u}{r} \cp \lens{f^\#}{f}, \tilde p \cp i'\right) \in \Rsm\left(\lens{I'}{O'}, M' \right)$$ where $S\times_M Q$ is the pullback
    \begin{center}
\begin{tikzcd}
                S \times_M Q \arrow[r, "\tilde p", dashed] \arrow[d, "\tilde i ", dashed,  swap]\arrow[dr, phantom, "\lrcorner", very near start]
                    & Q  \arrow[d, "i"]\\
                S \arrow[r, "p", swap]
                    & M\\
            \end{tikzcd}
\end{center}
and the induced lens is
\begin{center}
    \begin{tikzcd}
        \lensbig{T(S \times_M Q)}{S\times_M Q} \arrow[r, "\tilde i "', shift right=2]
        & \lensbig{TS}{S} \arrow[l, "T(\tilde i )^\#"', shift right = 2] \arrow[r, "r "', shift right=2]
        & \lensbig{I}{O} \arrow[r, "f "', shift right=2]  \arrow[l, "u"', shift right = 2]
        & \lensbig{I'}{O'}.\arrow[l, "f^\#"', shift right = 2]
    \end{tikzcd}
\end{center}
\end{itemize}
\end{thm}

\begin{proof}
The functor $T: \Ca \to \Lens_{\Aa}$ defines a profunctor $$\Lens_{\Aa}(T(-), -): \Lens_{\Aa} \topro \Ca$$
and the inclusion $J: \Ca \to \Span_{\Ca}$ defines a profunctor $$\Span_{\Ca}(-,J(-)): \Ca \topro \Span_{\Ca}.$$
Let $\Rsm$ be the composition of profunctors $$\Lens_{\Aa}(T(-), -) \cp \Span_{\Ca}(-,J(-)).$$ Lemmas~\ref{lem:rsm-objects} and~\ref{lem:rsm-morphs} show that $\Rsm$ has the desired behavior on objects and morphisms respectively.

Lastly, we want to define a laxator for $\Rsm$. Let $\phi$ be the op-laxator for $T$. Let $$\left(\lens{I}{O}, M\right), \left(\lens{I'}{O'}, M'\right) \in \ob \Lens_{\Aa} \times \Span_{\Ca}.$$ Define the laxator of $\Rsm$ so that the pair $$\left( \left(S, \lens{u}{r}, p\right), \left(S', \lens{u'}{r'}, p'\right)\right) \in \Rsm \left(\lens{I}{O}, M\right) \times \Rsm \left(\lens{I'}{O'}, M'\right)$$ maps to $$\left(S \times S', \phi_{S, S'} \cp \left( \lens{u}{r} \otimes \lens{u'}{r'}\right), p \times p'\right) \in \Rsm\left( \lens{I}{O} \otimes \lens{I'}{O'} , M \times M'\right).$$

The unitor is defined by the choice $$( 1_{\Ca}, \eta, \id_{1_{\Ca}}) \in \Rsm(1_{\Lens_{\Aa}}, 1_{\Ca})$$ where $\eta$ is the op-unitor for $T$.

The associativity and unitality conditions follow straightforwardly from the op-associativity and op-unitality conditions for $T$.
\end{proof}

The domain $\Lens_{\Aa} \times \Span_{\Ca}$ of the functor $\Rsm$ defines a syntax for composing open dynamical systems. Loosely, we depict objects $\left(\lens{I}{O}, M\right)$ as boxes with input wires corresponding to $I$, output wires corresponding to $O$, and exposed wires corresponding to $M$. See Figure~\ref{fig:rsm-boxes}.

A box of this type is filled with an element $\left( S, \lens{u}{r}, p\right)  \in \Rsm\left(\lens{I}{O}, M\right)$ which is interpreted as
\begin{itemize}
    \item $S$ --- a state space or set of states
    \item $\lens{u}{r} : \lens{TS}{S} \leftrightarrows \lens{I}{O}$ ---  dynamics  and a read-out function
    \item $p: S\to M$ --- the observation of the states from the perspective of the exposed ports
\end{itemize}

\begin{figure}\label{fig:rsm-boxes}
    \centering

    \subfloat[]{\begin{tikzpicture}[every text node part/.style={align=center}, oriented WD, bbx = 1cm, bby =.5cm, bb min width=2cm, bb min height=2cm,bb port length=4pt, bb port sep=1]

\node[bb={0}{3}{2}{4}] (B){};

\node[above = 4pt of B] () {$I$};
\node[below = 4pt of B] () {$O$};
\node[right = 4pt of B] () {$M$};

\end{tikzpicture}}\qquad
    \subfloat[]{\begin{tikzpicture}[every text node part/.style={align=center}, oriented WD, bbx = 1cm, bby =.5cm, bb min width=3cm, bb min height=3cm,bb port length=4pt, bb port sep=1]

\node[bb={0}{3}{2}{4}] (B){$S \in \ob \Ca$\\$\lens{u}{r}: \lens{TS}{S} \leftrightarrows \lens{I}{O}$\\$p: S \to M$};

\node[above = 4pt of B] () {$I$};
\node[below = 4pt of B] () {$O$};
\node[right = 4pt of B] () {$M$};

\end{tikzpicture}}

    \caption{(a) Loosely, this box represents the type $\left({I \choose O}, M\right) \in \Lens_{\Aa} \times \Span_{\Ca}.$ \newline (b) Such a box may be filled with an element\newline $\left(S, {u \choose r}, p\right) \in \Rsm\left({I \choose O}, M\right)$ \newline representing a choice of states, dynamics, read-out, and observation from exposed ports. }
\end{figure}

For a morphism $$\left(\lens{f^\#}{f}, M \from Q \to M'\right): \left(\lens{I}{O}, M\right) \to \left(\lens{I'}{O'}, M'\right)$$ in $\Lens_{\Aa} \times \Span_{\Ca}$, the set map $\Rsm\left(\lens{f^\#}{f}, M \from Q \to M'\right)$ defines the effects of the machine-style composition given by the lens $\lens{f^\#}{f}$ and resource sharing given by the span $M \from Q \to M'$ to dynamical systems in $\Rsm\left(\lens{I}{O}, M\right)$.

\newthought{In Section~\ref{sec:rsm-examples} we constructed}  the Lokta-Volterra  predator-prey model as the composition of four simple continuous dynamical systems. To see Theorem~\ref{thm:rsm} in action we will formalize this construction in the language of the algebra $$\Rsm: \Lens_{\Riem_{\Sub}/(-)} \times \Span_{\Riem} \to \Set$$ induced by the contravariant dynamical system doctrine consisting of the  indexed category $$\Riem_{\Sub}/(-): \Riem \op \to \Cat$$ described in Example~\ref{ex:rsm-cds-generalized-lens} and the section $T: \Riem \to \Lens_{\Riem_{\Sub}/(-)}$ described in Example~\ref{ex:rsm-cds-section-T}.

\begin{figure}\label{fig:rsm-CDS-total}
    \centering
    \begin{tikzpicture}[every text node part/.style={align=center}, oriented WD, bbx = 1cm, bby =.5cm, bb min width=1cm, bb min height=1cm,bb port length=4pt, bb port sep=1]

\node[bb={1}{0}{0}{1}] (rg) {};
\node[bb={0}{1}{1}{0}, right = of rg] (fg) {};
\node[bb={1}{0}{1}{0}, right= of fg] (rd) {};
\node[bb={0}{1}{0}{1}, right = of rd] (fd) {};

\node[bb={0}{0}{0}{0}, fit = {($(rg.north west) + (0, 1)$)($(fd.south east) - (0, 1)$)}] (tot) {};

\draw[ar, color =violet] let \p1=(rg.south east), \p2=(fg.north west), \n1=\bbportlen, \n2=\bby in
  	(rg_bot1') to [out=-90, in=-90](\x1/2+\x2/2,\y1-\n1) -- (\x1/2+\x2/2,\y2+\n1) to [out=90,in=90](fg_top1);
  	
\draw[ar, color =violet] let \p1=(fd.south east), \p2=(rd.north west), \n1=\bbportlen, \n2=\bby in
  	(fd_bot1') to [out=-90, in=-90](\x1/2+\x2/2,\y1-\n1) -- (\x1/2+\x2/2,\y2+\n1) to [out=90,in=90](rd_top1);
  	
\draw[color=cyan] let \p1=(rg.south west),  \p2=(fg.south east), \n1=\bbportlen, \n2=\bby in
    (rg_in1) to [out = 180, in = 180] (\x1 - \n1, \y1 - \n2-\n1) -- (\x2 + \n1, \y2 - \n2-\n1) to [out=0, in = 180] (rd_in1);
    
\draw[color=cyan] let \p1=(fd.north east),  \p2=(rd.north west), \n1=\bbportlen, \n2=\bby in
    (fd_out1) to [out = 0, in = 0] (\x1 + \n1, \y1 + \n2+\n1) -- (\x2 - \n1, \y2 + \n2+\n1) to [out=180, in = 0] (fg_out1);

\end{tikzpicture}
    \caption{The syntax for composing four continuous dynamical systems (represented by the inner boxes) into a single total system (represented by the outer box). The total system has trivial inputs, outputs, and exposed ports. }
\end{figure}

Figure~\ref{fig:rsm-CDS-total} defines a syntax for composing four continuous dynamical systems. A box with $I$ input ports, $O$ output ports, and $M$ exposed ports has type $$\left(\lens{\Rb^I \times \Rb^O \xrightarrow{\pi_2} \Rb^O}{\Rb^O}, \Rb^M\right) \in \ob\left( \Lens_{\Riem_{\Sub}/(-)} \times \Span_{\Riem} \right).$$ In the following we abuse notation by suppressing the submersion in each object of $\Lens_{\Riem_{\Sub}/(-)}$ since they are all given by projection onto the second coordinate, and instead we write $$\left(\lens{\Rb^I \times \Rb^O}{\Rb^O}, \Rb^M\right) \in \ob\left( \Lens_{\Riem_{\Sub}/(-)} \times \Span_{\Riem} \right).$$

The left-most and right-most inner box have $0$ input ports, $1$ output port, and $1$ exposed port. Therefore, these boxes have type $\left(\lens{\{*\} \times \Rb}{\Rb}, \Rb\right)$. On the other hand, the middle inner boxes have type $\left(\lens{\Rb \times \{*\} }{\{*\}}, \Rb\right)$. The outer box has trivial type $\left(\lens{\{*\} \times \{*\} }{\{*\}}, \{*\}\right)$.

The wiring in Figure~\ref{fig:rsm-CDS-total} represents a morphism
\begin{center}
    $\left(\lens{\{*\} \times \Rb}{\Rb}, \Rb\right)  \otimes \left(\lens{\Rb \times \{*\}}{\{*\}}, \Rb\right) \otimes \left(\lens{\Rb \times \{*\}}{\{*\}}, \Rb\right) \otimes \left(\lens{\{*\} \times \Rb}{\Rb}, \Rb\right) \to \left(\lens{\{*\} \times \{*\}}{\{*\}}, \{*\}\right)$
\end{center}in $\Lens_{\Riem_{\Sub}/(-)} \times \Span_{\Riem}$. Unwinding definitions, the wiring is given by the pair of morphisms:
\begin{enumerate}
    \item  the lens $\lens{\id \times \id}{!}: \lens{\Rb^2 \times \Rb^2}{\Rb^2} \leftrightarrows \lens{\{*\}}{\{*\}}$. The morphisms $\phi^{\mathsf{in}}=\id \times \id$ and $\phi^{\mathsf{out}} = !$ fit into the  commutative diagram on below.
    \begin{center}
    \begin{tikzcd}
        & \Rb^2 \arrow[r, dashed] \arrow[d, dashed] \arrow[dl, swap, "\phi^{\mathsf{in}} = \id \times \id", color = violet]\arrow[dr, phantom, "\lrcorner", very near start]
        & \{*\}\arrow[d, "\id"]\\
    \Rb^2 \times \Rb^2 \arrow[r, "\pi_2", swap]
        & \Rb^2 \arrow[r, "\phi^{\mathsf{out}} = !", swap, color = orange]
        & \{*\}
    \end{tikzcd}
    \end{center}
    \item the span {\color{cyan}$\Rb^4 \xleftarrow{\Delta} \Rb^2 \xrightarrow{!} \{*\}$}.
\end{enumerate}

Now that the syntax is established we can interpret each domain type\sidenote{Recall that a domain type is represented by an inner box in Figure~\ref{fig:rsm-CDS-total}.} as a continuous dynamical system. Then  composition along the wiring defined above results in an interpretation for the codomain type.\sidenote{The codomain type is represent by the outer box in Figure~\ref{fig:rsm-CDS-total}}

The left-most box in Figure~\ref{fig:rsm-CDS-total} has type $\left(\lens{\{*\} \times \Rb}{\Rb}, \Rb\right)$ and hence may be filled with an element $\left(S, \lens{u}{r}, p\right) \in \Rsm\left(\lens{\{*\} \times \Rb}{\Rb}, \Rb\right)$. The choice $$\left(\Rb, \lens{u}{\id}: \lens{T\Rb}{\Rb} \leftrightarrows \lens{\{*\} \times \Rb}{\Rb}, \id: \Rb \to \Rb\right)$$ where $u(r) = \beta r \in T_r\Rb$ models rabbit population growth. As short-hand, we visualize this filling as:
\vspace{-5pt}
\begin{figure}
    \centering
    \begin{tikzpicture}[every text node part/.style={align=center}, oriented WD, bbx = 1cm, bby =.5cm, bb min width=1cm, bb min height=1cm,bb port length=4pt, bb port sep=1]

\node[bb={1}{0}{0}{1}] (rg) {$\dot r= \beta r$};
\node[above = 1pt of rg_in1] () {\scriptsize $r$};
\node[left = 1pt of rg_bot1', yshift=-1pt] () {\scriptsize $r$};

\end{tikzpicture}
\end{figure}
\vspace{-10pt}

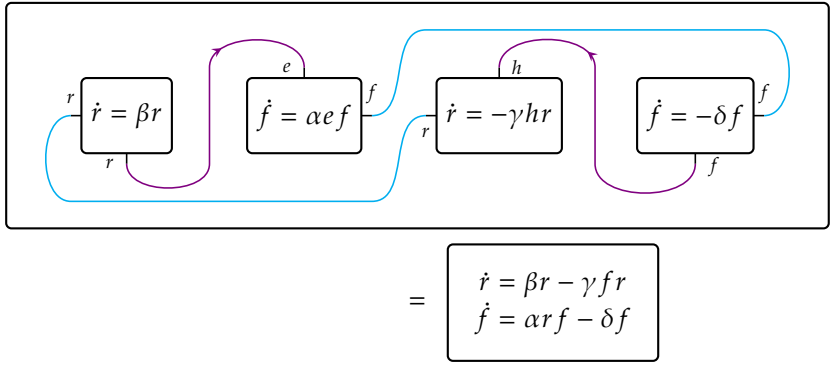
\begin{figure}\label{fig:rsm-CDS-total-filled}
    \centering
    \begin{tikzpicture}[every text node part/.style={align=center}, oriented WD, bbx = 1cm, bby =.5cm, bb min width=1cm, bb min height=1cm,bb port length=4pt, bb port sep=1]

\node[bb={1}{0}{0}{1}] (rg) {$\dot r= \beta r$};
\node[bb={0}{1}{1}{0}, right = of rg] (fg) {$\dot f= \alpha e f$};
\node[bb={1}{0}{1}{0}, right= of fg] (rd) {$\dot r = -\gamma hr$};
\node[bb={0}{1}{0}{1}, right = of rd] (fd) {$\dot f = -\delta f$};

\node[bb={0}{0}{0}{0}, fit = {($(rg.north west) + (0, 1)$)($(fd.south east) - (0, 1)$)}] (tot) {};

\draw[ar, color =violet] let \p1=(rg.south east), \p2=(fg.north west), \n1=\bbportlen, \n2=\bby in
  	(rg_bot1') to [out=-90, in=-90](\x1/2+\x2/2,\y1-\n1) -- (\x1/2+\x2/2,\y2+\n1) to [out=90,in=90](fg_top1);
  	
\draw[ar, color =violet] let \p1=(fd.south east), \p2=(rd.north west), \n1=\bbportlen, \n2=\bby in
  	(fd_bot1') to [out=-90, in=-90](\x1/2+\x2/2,\y1-\n1) -- (\x1/2+\x2/2,\y2+\n1) to [out=90,in=90](rd_top1);
  	
\draw[color=cyan] let \p1=(rg.south west),  \p2=(fg.south east), \n1=\bbportlen, \n2=\bby in
    (rg_in1) to [out = 180, in = 180] (\x1 - \n1, \y1 - \n2-\n1) -- (\x2 + \n1, \y2 - \n2-\n1) to [out=0, in = 180] (rd_in1);
    
\draw[color=cyan] let \p1=(fd.north east),  \p2=(rd.north west), \n1=\bbportlen, \n2=\bby in
    (fd_out1) to [out = 0, in = 0] (\x1 + \n1, \y1 + \n2+\n1) -- (\x2 - \n1, \y2 + \n2+\n1) to [out=180, in = 0] (fg_out1);
    
\node[above = 1pt of rg_in1] () {\scriptsize $r$};
\node[left = 1pt of rg_bot1'] () {\scriptsize $r$};

\node[left = 1pt of fg_top1] () {\scriptsize $e$};
\node[above = 1pt of fg_out1] () {\scriptsize $f$};

\node[right = 1pt of rd_top1] () {\scriptsize $h$};
\node[below = 1pt of rd_in1] () {\scriptsize $r$};

\node[right = 1pt of fd_bot1', yshift = -1pt] () {\scriptsize $f$};
\node[above = 1pt of fd_out1] () {\scriptsize $f$};

\node[below = 5ex of tot] (eq) {$=$};
\node[right = 1ex of eq, bb={0}{0}{0}{0},inner sep = 10pt] () {$\dot r = \beta r - \gamma fr$\\$\dot f =\alpha rf - \delta f$};

\end{tikzpicture}
    \caption{On the left-hand side of the equation we fill the syntax defined in  Figure~\ref{fig:rsm-CDS-total} with continuous dynamical systems which from left to right correspond to rabbit population growth, fox population growth, rabbit population decline, and fox population decline. On the right-hand side of the equation, the resulting interpretation induced by $\Rsm$ is the Lokta-Volterra predator-prey model.}
\end{figure}

Likewise, we may fill the remaining boxes with open continuous dynamical systems that model fox population growth, rabbit population decline,  and fox population decline. The short-hand for these systems fill the inner boxes in the left-hand side of the equation in  Figure~\ref{fig:rsm-CDS-total-filled}. The set map $\Rsm\left(\lens{\id \times \id}{!}, \Rb^4 \xleftarrow{\Delta} \Rb^2 \xrightarrow{!} \{*\}\right)$\sidenote{Recall that the wirings in Figures~\ref{fig:rsm-CDS-total} and~\ref{fig:rsm-CDS-total-filled} represent the morphism $\left(\lens{\id \times \id}{!}, \Rb^4 \xleftarrow{\Delta} \Rb^2 \xrightarrow{!} \{*\}\right)$.} maps the quadruple of fillings to the continuous dynamical system
$$\dot r = \beta r - \gamma f r, \quad \dot f = \alpha r f - \delta f$$ the Lokta-Volterra predator-prey model.

\newpage
\section{Proof of Main Theorem}
Throughout let 
\begin{itemize}
    \item $\Ca$ be a cartesian monoidal category with pullbacks
    \item $\Aa: \Ca\op \to \Cat$ be a monoidal indexed category
    \item $T: \Ca \to \Lens_{\Aa}$  be an oplax section of the forgetful functor with colaxator $\phi$
\end{itemize}

\begin{lem}\label{lem:rsm-objects}
Let $\lens{I}{O} \in \ob \Lens_{\Aa}$ and $M \in \ob \Span_{\Ca}$. Define $\Rsm\left(\lens{I}{O}, M\right)$ to be the set of triples $\left( S, \lens{u}{r}, p \right)$ with $S \in \ob \Ca$, $\lens{u}{r}: \lens{TS}{S} \leftrightarrows \lens{I}{O}$ in $\Lens_{\Aa}$, and $p: S \to M$ in $\Ca$. Then, $\Rsm\left(\lens{I}{O}, M\right)$ is isomorphic to the coend $$\int^{S \in \Ca} \Lens_{\Aa}\left(\lens{TS}{S}, \lens{I}{O}\right) \times \Span_{\Ca}(M, S).$$
\end{lem}

\begin{proof}
For each $S \in \ob \Ca$ define a set map 
$$\omega_S: \Lens_{\Aa}\left(\lens{TS}{S}, \lens{I}{O}\right) \times \Span_{\Ca}(M, S) \to \Rsm\left(\lens{I}{O}, M\right)$$ 
by $$\omega_S\left(\lens{u}{r}, M \xleftarrow{p} S' \xrightarrow{q} S\right) = \left( S', \lens{Tq^\#}{q} \cp \lens{u}{r}, p\right).$$ 
First we will show that $$\omega: \Lens_{\Aa}\left(T(-), \lens{I}{O}\right) \times \Span_{\Ca}(M, -) \Rightarrow \Rsm\left(\lens{I}{O}, M\right)$$
is a cowedge of $$\Lens_{\Aa}\left(T(-), \lens{I}{O}\right) \times \Span_{\Ca}(M, -): \Ca\op \times \Ca \to \Set.$$

 Let $f: S_1\to S_2$ in $\Ca$. The following diagram commutes by unwinding definitions.
\begin{center}
\begin{tikzcd}
    \Lens_{\Aa}\left(\lens{TS_2}{S_2}, \lens{I}{O}\right) \times \Span_{\Ca}(M, S_1) \arrow[r] \arrow[d]
    
    & \Lens_{\Aa}\left(\lens{TS_1}{S_1}, \lens{I}{O}\right) \times \Span_{\Ca}(M, S_1) \arrow[d, "\omega_{S_1}"]\\
    \Lens_{\Aa}\left(\lens{TS_2}{S_2}, \lens{I}{O}\right) \times \Span_{\Ca}(M, S_2) \arrow[r, "\omega_{S_2}", swap]
    & \Rsm\left(\lens{I}{O}, M\right)\\
\end{tikzcd}
\end{center}

Next we want to show that $\Rsm\left(\lens{I}{O}, M\right)$ is the universal cowedge. Suppose $\alpha: \Lens_{\Aa}\left(T(-), \lens{I}{O}\right) \times \Span_{\Ca}(M, -) \Rightarrow Y$ is any cowedge. We want to show there exists $h: \Rsm\left(\lens{I}{O}, M\right)\to Y$ such that for all $S \in \Ca$, the triangle below commutes.
\begin{center}
\begin{tikzcd}
    \Lens_{\Aa}\left(\lens{TS}{S}, \lens{I}{O}\right) \times \Span_{\Ca}(M, S) \arrow[d, swap, "\omega_S"] \arrow[r, "\alpha_S"] 
    & Y \\
    \Rsm\left(\lens{I}{O}, M\right)\arrow[ur, "h", swap, bend right = 20, dashed]
    & \\
\end{tikzcd}
\end{center}

For $\left( S, \lens{u}{r}, p\right) \in \Rsm\left(\lens{I}{O}, M\right)$ define 
\begin{equation}\label{eq:rsm-universalprop}
    h\left( S, \lens{u}{r}, p \right) = \alpha_{S}\left( \lens{u}{r}, M \xleftarrow{p} S \xrightarrow{\id} S\right)
\end{equation}

Then for all $$\left( \lens{u}{r}, M \xleftarrow{p} S' \xrightarrow{q} S \right) \in  \Lens_{\Aa}\left(\lens{TS}{S}, \lens{I}{O}\right) \times \Span_{\Ca}(M, S),$$ we have
\begin{align*}
    (h \circ \omega_S)\left( \lens{u}{r}, M \xleftarrow{p} S' \xrightarrow{q} S \right)
    & = h\left( S', \lens{Tq^\#}{q} \cp \lens{u}{r}, p\right)\\
    & = \alpha_{S'}\left(\lens{Tq^\#}{q} \cp \lens{u}{r}, M \xleftarrow{p} S' \xrightarrow{\id} S'\right)\\
    & = \alpha_{S}\left( \lens{u}{r}, M \xleftarrow{p} S' \xrightarrow{q} S \right)
\end{align*} where the last line follows from the commuting square   induced by $q: S' \to S$ below.
\begin{center}
\begin{tikzcd}
    \Lens_{\Aa}\left(\lens{TS}{S}, \lens{I}{O}\right) \times \Span_{\Ca}(M, S') \arrow[r] \arrow[d]
    
    & \Lens_{\Aa}\left(\lens{TS'}{S'}, \lens{I}{O}\right) \times \Span_{\Ca}(M, S') \arrow[d, "\alpha_{S'}"]\\
    \Lens_{\Aa}\left(\lens{TS}{S}, \lens{I}{O}\right) \times \Span_{\Ca}(M, S) \arrow[r, "\alpha_{S}", swap]
    & \Rsm\left(\lens{I}{O}, M\right)\\
\end{tikzcd}
\end{center}

Lastly we must show that $h$ is unique. Suppose $\tilde h : \Rsm\left(\lens{I}{O}, M\right)\to Y$ also satisfies $\tilde h \circ \omega_S = \alpha_S$ for all $S \in \ob \Ca$. For $\left( S, \lens{u}{r}, p\right) \in \Rsm\left(\lens{I}{O}, M\right)$,  $$\left( S, \lens{u}{r}, p\right) = \omega_S\left(\lens{u}{r},M \xleftarrow{p} S \xrightarrow{\id} S\right)$$ implies $$\tilde h\left( S, \lens{u}{r}, p\right) = \alpha_S\left(\lens{u}{r},M \xleftarrow{p} S \xrightarrow{\id} S\right) = h\left( S, \lens{u}{r}, p \right).$$
\end{proof}

\begin{lem}\label{lem:rsm-morphs}
Let $\lens{f^\#}{f}: \lens{I}{O} \leftrightarrows \lens{I'}{O'}$ be a morphism in $\Lens_{\Aa}$ and let $M \xleftarrow{i} Q \xrightarrow{i'} M'$ be a span in $\Ca$. These induce a natural transformation
\begin{align*}
    \Lens_{\Aa}\left(T(-), \lens{f^\#}{f}\right) \times \Span_{\Ca}( M \xleftarrow{i} Q \xrightarrow{i'} M', -) : &\Lens_{\Aa}\left(T(-), \lens{I}{O}\right) \times \Span_{\Ca}(M, -) \\
     \Rightarrow & \Lens_{\Aa}\left(T(-), \lens{I'}{O'}\right) \times \Span_{\Ca}(M', -) .
\end{align*}

Let  $\Rsm\left(\lens{f^\#}{f}, M \xleftarrow{i} Q \xrightarrow{i'} M' \right)$ be the set map which takes the triple $(S, \lens{u}{r}, p) \in \Rsm\left(\lens{I}{O}, M\right)$ to the triple $$\left(S \times_M Q, \lens{T\tilde i ^\#}{\tilde i } \cp \lens{u}{r} \cp \lens{f^\#}{f}, \tilde p \cp i'\right) \in \Rsm\left(\lens{I'}{O'}, M' \right).$$ 

Then, $\Rsm\left(\lens{f^\#}{f}, M \xleftarrow{i} Q \xrightarrow{i'} M' \right)$ is isomorphic to  $$\int^{S \in \Ca}\Lens_{\Aa}\left(T(-), \lens{f^\#}{f}\right) \times \Span_{\Ca}( M \xleftarrow{i} Q \xrightarrow{i'} M', -).$$
\end{lem}
\begin{proof}
The transformation $$\Lens_{\Aa}\left(T-, \lens{f^\#}{f}\right) \times \Span_{\Ca}( M \xleftarrow{i} Q \xrightarrow{i'} M', -)$$
induces a cowedge $$\alpha: \Lens_{\Aa}\left(T-, \lens{I}{O}\right) \times \Span_{\Ca}(M, -) \Rightarrow \Rsm\left(\lens{I'}{O'}, M'\right)$$ defined such that $$\alpha_S: \Lens_{\Aa}\left(\lens{TS}{S}, \lens{I}{O}\right) \times \Span_{\Ca}(M, S) \Rightarrow \Rsm\left(\lens{I'}{O'}, M'\right)$$ maps the pair $\left(\lens{u}{r}: \lens{TS}{S} \leftrightarrows\lens{I}{O}, M \xleftarrow{p} S' \xrightarrow{q} S\right)$ to the triple $$\left( S' \times_M Q, \lens{T(S' \times_M Q)}{S' \times_M Q} \lensmap{T(\tilde{i} \cp p)^\#}{\tilde i  \cp p} \lens{TS}{S} \lensmap{u}{r} \lens{I}{O} \lensmap{f^\#}{f} \lens{I'}{O'}, S' \times_M Q \to Q \xrightarrow{i'} M' \right).$$

Then $$\int^{S \in \Ca}\Lens_{\Aa}\left(T(-), \lens{f^\#}{f}\right) \times \Span_{\Ca}( M \xleftarrow{i} Q \xrightarrow{i'} M', -)$$ is the unique map $$\Rsm\left(\lens{I}{O}, M\right) \to\Rsm\left(\lens{I'}{O'}, M'\right)$$ such that the diagram below commutes for all $S \in \Ca$, 

\begin{center}
\begin{tikzcd}
    \Lens_{\Aa}\left(\lens{TS}{S}, \lens{I}{O}\right) \times \Span_{\Ca}(M, S) \arrow[d, swap, "\omega_S"] \arrow[r, "\alpha_S"] 
    & \Rsm\left(\lens{I'}{O'}, M'\right) \\
    \Rsm\left(\lens{I}{O}, M\right)\arrow[ur, dashed, swap, bend right = 20]
    & \\
\end{tikzcd}
\end{center}
Following Equation~\ref{eq:rsm-universalprop} and unwinding definitions, $$\int^{S \in \Ca}\Lens_{\Aa}\left(T(-), \lens{f^\#}{f}\right) \times \Span_{\Ca}( M \xleftarrow{i} Q \xrightarrow{i'} M', -)$$ takes the triple $\left(S, \lens{u}{r}, p\right) \in \Rsm\left(\lens{I}{O}, M\right)$ to 
$$\alpha_S\left(\lens{u}{r}, M \xleftarrow{p} S \xrightarrow{\id} S\right)=  \Rsm\left(\lens{f^\#}{f}, M \xleftarrow{i} Q \xrightarrow{i'} M' \right)\left(S, \lens{u}{r}, p\right) \in \Rsm\left(\lens{I'}{O'}, M'\right)$$
\end{proof}

\bibliography{ref}

\begin{thebibliography}{12}
\providecommand{\natexlab}[1]{#1}
\providecommand{\url}[1]{\texttt{#1}}
\expandafter\ifx\csname urlstyle\endcsname\relax
  \providecommand{\doi}[1]{doi: #1}\else
  \providecommand{\doi}{doi: \begingroup \urlstyle{rm}\Url}\fi

\bibitem[Baez and Pollard(2017)]{Baez_2017}
John~C. Baez and Blake~S. Pollard.
\newblock A compositional framework for reaction networks.
\newblock \emph{Reviews in Mathematical Physics}, 29\penalty0 (09), 2017.

\bibitem[Baez et~al.(2016)Baez, Fong, and Pollard]{baez2016compositional}
John~C Baez, Brendan Fong, and Blake~S Pollard.
\newblock A compositional framework for markov processes.
\newblock \emph{Journal of Mathematical Physics}, 57\penalty0 (3):\penalty0
  033301, 2016.

\bibitem[Fong and Spivak(2018)]{fong2018hypergraph}
Brendan Fong and David~I Spivak.
\newblock Hypergraph categories, 2018.

\bibitem[Fong and Spivak(2019)]{fong2019invitation}
Brendan Fong and David~I Spivak.
\newblock \emph{An invitation to applied category theory: seven sketches in
  compositionality}.
\newblock Cambridge University Press, 2019.

\bibitem[Krohn and Rhodes(1965)]{krohn1965algebraic}
Kenneth Krohn and John Rhodes.
\newblock Algebraic theory of machines. i. prime decomposition theorem for
  finite semigroups and machines.
\newblock \emph{Transactions of the American Mathematical Society},
  116:\penalty0 450--464, 1965.

\bibitem[Leinster(2004)]{leinster2004higher}
Tom Leinster.
\newblock \emph{Higher operads, higher categories}, volume 298.
\newblock Cambridge University Press, 2004.

\bibitem[Moeller and Vasilakopoulou(2018)]{moeller2018monoidal}
Joe Moeller and Christina Vasilakopoulou.
\newblock Monoidal grothendieck construction, 2018.

\bibitem[Myers(2020)]{myers2020double}
David~Jaz Myers.
\newblock Double categories of open dynamical systems (extended abstract),
  2020.

\bibitem[Schultz et~al.(2016)Schultz, Spivak, and
  Vasilakopoulou]{schultz2016dynamical}
Patrick Schultz, David~I. Spivak, and Christina Vasilakopoulou.
\newblock Dynamical systems and sheaves, 2016.

\bibitem[Spivak(2019)]{spivak2019generalized}
David~I. Spivak.
\newblock Generalized lens categories via functors $\mathcal{C}^{\rm
  op}\to\mathsf{Cat}$, 2019.

\bibitem[Spivak(2020)]{spivak2020poly}
David~I Spivak.
\newblock Poly: An abundant categorical setting for mode-dependent dynamics.
\newblock \emph{arXiv preprint arXiv:2005.01894}, 2020.

\bibitem[Vagner et~al.(2014)Vagner, Spivak, and Lerman]{vagner2014algebras}
Dmitry Vagner, David~I Spivak, and Eugene Lerman.
\newblock Algebras of open dynamical systems on the operad of wiring diagrams.
\newblock \emph{arXiv preprint arXiv:1408.1598}, 2014.

\end{thebibliography}
\bibliographystyle{plainnat}

\end{document}